\documentclass[11pt]{amsart}
\usepackage[T1]{fontenc}
\usepackage{lmodern}
\usepackage{amsmath,amssymb,amsthm,booktabs}
\usepackage[margin=1in]{geometry}
\usepackage[colorlinks=true,linkcolor=blue,citecolor=blue,urlcolor=blue]{hyperref}
\hypersetup{pdftitle={Complete Bernstein functions and real powers of P-matrices},pdfauthor={Wei Xie},pdfkeywords={P-matrix, complete Bernstein function, principal matrix power, total positivity}}
\newtheorem{theorem}{Theorem}[section]
\newtheorem{proposition}[theorem]{Proposition}
\newtheorem{lemma}[theorem]{Lemma}
\newtheorem{corollary}[theorem]{Corollary}
\newcommand{\R}{\mathbb R}
\newcommand{\Pn}{\mathcal P_n}
\newcommand{\TP}{\mathrm{TP}}
\newcommand{\En}{\mathcal E_n}
\newcommand{\diag}{\operatorname{diag}}
\newcommand{\tr}{\operatorname{tr}}
\newcommand{\abs}[1]{\left|#1\right|}
\newcommand{\norm}[1]{\left\|#1\right\|}
\numberwithin{equation}{section}
\title{Complete Bernstein functions and real powers of P-matrices}
\author{Wei Xie}
\date{October 25, 2026}

\subjclass[2020]{15A16, 15B48}
\keywords{P-matrix, complete Bernstein function, principal matrix power, total positivity}

\begin{document}
\begin{abstract}
We classify the nonzero complete Bernstein functions whose matrix functions
preserve every real P-matrix in every order. They are precisely
$f(z)=a+bz$, with $a,b\ge0$ and $a+b>0$, and
$f(z)=a+cz/(z+s)$, with $a\ge0$ and $c,s>0$.
For each fixed order $n\ge4$, the preservers vanishing at zero are exactly
$cz$ and $cz/(z+s)$; every other such function admits a strictly totally
positive counterexample. For a finite positive measure $\mu$ on $(0,\infty)$, put
$g(z)=\int z/(z+s)\,\mu(ds)$ and $M_t=\int s^t\,\mu(ds)$.
In every fixed order $n\ge4$, a finite constant shift restores
preservation on all real P-matrices if and only if $M_{r_n}$ and
$M_{-r_n}$ are finite, where $r_n=\lfloor n/2\rfloor-1$.
Every $a>c_n\sqrt{M_{r_n}M_{-r_n}}$ is sufficient, with $c_n$
depending only on $n$.
A complete Bernstein function vanishing at zero admits a possibly different repairing
constant in every order exactly when it is linear, or is a pure kernel
mixture with all positive and negative integer moments finite.
For principal real powers, the preserving exponents are $\mathbb R$
in order one, $[-1,1]$ in orders two and three, and $\{-1,0,1\}$
in every higher order. An explicit rational determinant-one totally
positive family gives negative principal minors uniformly on prescribed
compact subsets of $(-1,0)\cup(0,1)$.
\end{abstract}
\maketitle

\section{Introduction}

A real square matrix is a P-matrix if all its nonempty principal minors
are positive. We write $\Pn$ for the class of P-matrices of order $n$.
Real symmetric positive definite matrices remain positive definite under
every real power. Inversion also preserves P-matrices, by Jacobi's
complementary-minor identity.

A P-matrix has no eigenvalues in $(-\infty,0]$. For a complex number
$\lambda\notin(-\infty,0]$ and a real number $t$, its principal power is
\[
 \lambda^t=|\lambda|^t e^{it\arg\lambda},
 \qquad -\pi<\arg\lambda<\pi.
\]
If $A\in\Pn$ is diagonalizable, with eigenvalues
$\lambda_1,\ldots,\lambda_n$, choose a diagonalization
\[
 A=V\diag(\lambda_1,\ldots,\lambda_n)V^{-1}.
\]
Its principal real powers are
\[
 A^t=V\diag(\lambda_1^t,\ldots,\lambda_n^t)V^{-1}.
\]
The definition for matrices with nontrivial Jordan blocks is given in
Section~\ref{sec:low}.

Mondal, Sivakumar, and Tsatsomeros~\cite{MST} studied real powers of
P-matrices using Pick functions. Their Corollary~3.15 (p.~363) states that
\[
 A\in\Pn\quad\Longrightarrow\quad A^t\in\Pn
 \qquad(0<t\le1)
\]
in every order.
We revisit this assertion and ask, for each fixed order $n$, which real
exponents $t$ satisfy $A^t\in\Pn$ for every $A\in\Pn$.

For $s>0$, the rational function $z/(z+s)$ preserves P-matrices, since
\[
 A(A+sI)^{-1}=(I+sA^{-1})^{-1}
\]
and P-matrices are closed under inversion and addition of a nonnegative
scalar multiple of the identity; see also Fallat and
Tsatsomeros~\cite[Lemma~2.2 and Theorem~3.1]{FT}.
The standard integral representation of $z^t$ for $0<t<1$ expresses it
as a positive mixture of these rational functions. Such mixtures, together
with nonnegative constant and linear terms, form the class of complete
Bernstein functions. We consider which complete Bernstein functions
preserve P-matrices of a fixed order.

For a fixed complete Bernstein function $g$ with $g(0+)=0$ and a fixed
$A\in\Pn$, the matrix $aI+g(A)$ is a P-matrix for all sufficiently large
$a$, since $a^{-1}(aI+g(A))\to I$ as $a\to\infty$.
We ask under what conditions on $g$ there is a constant $a\ge0$ such that
\[
 aI+g(A)\in\Pn\qquad\text{for every }A\in\Pn.
\]

A matrix is called a PM-matrix in \cite{MST} if all its positive integer
powers are P-matrices. Theorem~4.1 of that paper (p.~363) states that every
principal real power of a PM-matrix is a P-matrix.
A matrix is strictly totally positive if all its nonempty minors are
positive; we denote this class by $\TP_n$.
The Cauchy--Binet formula gives
\[
 A\in\TP_n\quad\Longrightarrow\quad A^k\in\TP_n\subseteq\Pn
 \qquad(k=1,2,\ldots).
\]
Johnson and Walch~\cite[Theorem~3.8]{JW} proved that, for each order $n$,
there is a constant $\gamma_n$ such that $A^t$ is strictly totally positive
for every $A\in\TP_n$ and every $t>\gamma_n$.
For strictly totally positive inputs, we consider the P-matrix property
of fractional powers with $0<t<1$ and of the complete Bernstein functions
above.

\section{Preserving exponents and functions}\label{sec:low}

Put $[n]=\{1,\ldots,n\}$ and write $I^c=[n]\setminus I$ for
$I\subseteq[n]$. For increasing index sets $I,J$ of equal size,
$A[I,J]$ denotes the submatrix on those rows and columns. The empty
principal minor has value one.

For an invertible matrix $B$ of order $n$, Jacobi's complementary-minor
identity is
\begin{equation}\label{eq:jacobi}
 \det B^{-1}[I,J]
 =(-1)^{\sum_{i\in I}i+\sum_{j\in J}j}
   \frac{\det B[J^c,I^c]}{\det B}.
\end{equation}
In particular,
\begin{equation}\label{eq:jprincipal}
 \det B^{-1}[I,I]=\frac{\det B[I^c,I^c]}{\det B}.
\end{equation}

For $A\in\Pn$ and $s\ge0$, expansion in principal minors gives
\begin{equation}\label{eq:p}
p_A(s):=\det(sI+A)
=\sum_{J\subseteq[n]}\det A[J,J]s^{n-|J|}>0.
\end{equation}
This proves that $A$ has no nonpositive real eigenvalues.

For a Jordan block $J=\lambda I_m+N$, where $N^m=0$, define
\[
 J^t=\sum_{j=0}^{m-1}\binom tj\lambda^{t-j}N^j,\qquad
 \binom t0=1,\qquad
 \binom tj=\frac{t(t-1)\cdots(t-j+1)}{j!}\quad(j\ge1).
\]
Here each scalar power uses the branch specified in the Introduction.
If $A=S\diag(J_1,\ldots,J_q)S^{-1}$, define
\[
 A^t=S\diag(J_1^t,\ldots,J_q^t)S^{-1}.
\]
This is independent of the chosen Jordan representation: a polynomial
matching $z^t$ and its required derivatives at each eigenvalue gives the
same matrix. For real $A$, conjugate eigenvalues use conjugate scalar
powers, so $A^t$ is real. This definition agrees with ordinary integer
powers.

Define
\begin{equation}\label{eq:E}
\En=\{t\in\mathbb R:A^t\in\Pn\text{ for every }A\in\Pn\}.
\end{equation}

\begin{theorem}\label{thm:main}
The sets in \eqref{eq:E} are
\[
\mathcal E_1=\mathbb R,\qquad
\mathcal E_2=\mathcal E_3=[-1,1],\qquad
\mathcal E_n=\{-1,0,1\}\quad(n\ge4).
\]
For every fixed $0<|t|<1$, the following conditions are equivalent:
\begin{enumerate}
\item $A^t\in\Pn$ for every $A\in\Pn$;
\item $A^t\in\Pn$ for every $A\in\TP_n$;
\item $n\le3$.
\end{enumerate}
\end{theorem}

\begin{samepage}
\begin{theorem}\label{thm:uniform}
Let $n\ge4$ and $0<a\le b<1$, and put
$K=[-b,-a]\cup[a,b]$. There is a matrix
$A\in\TP_n\cap\mathrm{SL}_n(\mathbb R)$ such that
\begin{equation}\label{eq:negative}
\det A^t[\{1,4\},\{1,4\}]<0\qquad(t\in K).
\end{equation}
In order four, $A$ can be chosen with rational entries from the explicit
family in Section~\ref{sec:family}. For each fixed $K$, the matrices satisfying
\eqref{eq:negative} form a nonempty open subset of $\TP_n$ and a nonempty
relatively open subset after imposing determinant one.
\end{theorem}
\end{samepage}

For each fixed $A\in\Pn$, the matrix $A^t$ belongs to $\Pn$ for all
sufficiently small $|t|$, since $A^t\to I$ as $t\to0$ and $\Pn$ is open.

For $n\le3$, every nonempty principal minor has order $1$, $n-1$, or
$n$. In order four, the additional principal minors have order two.

\begin{proposition}\label{prop:low}
Let $A\in\Pn$ and $|t|\le1$. The principal minors of $A^t$ of orders
$1$, $n-1$, and $n$ are positive whenever these orders are nonzero.
Consequently $A^t\in\Pn$ if $n\le3$.
\end{proposition}

\begin{proof}
For $i\in[n]$, write $\widehat i=[n]\setminus\{i\}$.
The diagonal cofactor formula and \eqref{eq:p} imply
\begin{align}
((A+sI)^{-1})_{ii}
&=\frac{\det(sI+A[\widehat i,\widehat i])}{p_A(s)}>0,\label{eq:res1}\\
(A(A+sI)^{-1})_{ii}
&=1-s((A+sI)^{-1})_{ii}
=\frac{\displaystyle\sum_{J\ni i}\det A[J,J]s^{n-|J|}}{p_A(s)}>0.
\label{eq:res2}
\end{align}
For $0<t<1$, the scalar beta integral, analytically continued to
$\mathbb C\setminus(-\infty,0]$, gives
\begin{align}
A^t&=\frac{\sin\pi t}{\pi}
\int_0^\infty s^{t-1}A(A+sI)^{-1}\,ds,\label{eq:int1}\\
A^{-t}&=\frac{\sin\pi t}{\pi}
\int_0^\infty s^{-t}(A+sI)^{-1}\,ds.\label{eq:int2}
\end{align}
The scalar identities and their derivatives with respect to the scalar
argument give these formulas on each Jordan block.
The derivatives may be taken under the integrals: near
zero the resolvents and their derivatives are bounded, while at infinity
each derivative has at least the decay of the corresponding integrand.
The first integrand is $O(s^{t-2})$ and the second is $O(s^{-t-1})$ at
infinity, so both matrix integrals converge.
Equations~\eqref{eq:res1}--\eqref{eq:int2} show that the diagonal entries
of both $A^t$ and $A^{-t}$ are positive.

The eigenvalues of $A^u$, counted with multiplicity, are the principal
powers of the eigenvalues of $A$. A positive real eigenvalue $\lambda$
contributes $\lambda^u$ to the determinant, and a conjugate pair
$\lambda,\overline\lambda$ contributes $|\lambda|^{2u}$. Therefore
\begin{equation}\label{eq:detpower}
\det A^u=(\det A)^u>0\qquad(u\in\R).
\end{equation}
Taking $B=A^u$ in \eqref{eq:jprincipal}, with $u=t$ or $u=-t$, gives
\[
\det A^u[\widehat i,\widehat i]
=\det A^u\,(A^{-u})_{ii}>0.
\]
The endpoints $t=1$ and $t=0$ follow from the defining inequalities and
$A^0=I$. At $t=-1$, \eqref{eq:jprincipal} gives
$\det A^{-1}[J,J]=\det A[J^c,J^c]/\det A>0$ for every $J$.
If $n\le3$, the indicated orders exhaust the nonempty principal minors.
\end{proof}

For $|t|>1$, let $\theta=\pi/(|t|+1)$ and
\[
R_\theta=\begin{pmatrix}\cos\theta&-\sin\theta\\
\sin\theta&\cos\theta\end{pmatrix}.
\]
Since $0<\theta<\pi/2$, this is a P-matrix. Its principal power is
$R_{t\theta}$, whose diagonal entries equal
$\cos(t\theta)=\cos(\pi-\theta)=-\cos\theta<0$.
Taking the direct sum with $I_{n-2}$ gives a counterexample in every
order $n\ge2$.

\begin{corollary}\label{cor:TPouter}
If $A\in\TP_n$, then the principal minors of $A^t$ of orders
$1$, $n-1$, and $n$ are positive for every real $t$ whenever these
orders are nonzero. In particular, $A^t\in\Pn$ for every real $t$
if $n\le3$.
\end{corollary}

\begin{proof}
A strictly totally positive matrix has distinct positive eigenvalues.
Indeed, its $k$th compound, whose entries are its minors of order $k$,
is entrywise positive. Its eigenvalues are the products of $k$
eigenvalues of $A$. The Perron--Frobenius theorem, applied successively
to these compounds, gives $\lambda_1>\cdots>\lambda_n>0$.
Choose a positive integer $m\ge\max(1,|t|)$. Then $A^m\in\TP_n$,
and its positive spectrum gives $(A^m)^{t/m}=A^t$ on the principal
branch. Proposition~\ref{prop:low} applies to $A^m$ and $|t/m|\le1$.
\end{proof}

On $\mathbb C\setminus(-\infty,0]$, we also consider nonzero functions
of the form
\begin{equation}\label{eq:kernel}
 f(z)=a+bz+\int_{(0,\infty)}\frac{z}{z+s}\,\mu(ds),\qquad
 a,b\ge0,\quad \mu\ge0,\quad
 \int_{(0,\infty)}\frac{\mu(ds)}{1+s}<\infty.
\end{equation}
These are the nonzero complete Bernstein functions, with $f(0+)=a$;
see the standard representation in
Jacob and Schilling~\cite[formula~(13)]{JS}. For a real matrix $A$ whose spectrum avoids $(-\infty,0]$,
define
\begin{equation}\label{eq:matrixkernel}
 f(A)=aI+bA+\int_{(0,\infty)}A(A+sI)^{-1}\,\mu(ds).
\end{equation}
The matrix kernel is bounded near $s=0$ and is $O(s^{-1})$ at infinity,
so the integral converges. It is real and agrees with the ordinary
spectral function, including on nontrivial Jordan blocks: the scalar
kernel and each fixed derivative are bounded near zero and $O(s^{-1})$
at infinity, locally uniformly for $z\notin(-\infty,0]$, so the scalar
integral can be differentiated to evaluate each Jordan block.

\begin{theorem}\label{thm:functions}
Let $f$ be a nonzero function of the form \eqref{eq:kernel}.
For every $A\in\Pn$, the principal minors of $f(A)$ of orders
$1$, $n-1$, and $n$ are positive whenever these orders are nonzero.
Consequently $f(A)\in\Pn$ for $n\le3$.

Suppose in addition that $f(0+)=0$.
For each fixed $n\ge4$, the following conditions are equivalent:
\begin{enumerate}
\item $f(A)\in\Pn$ for every $A\in\Pn$;
\item $f(A)\in\Pn$ for every $A\in\TP_n$;
\item $f(z)=cz$ for some $c>0$, or
      $f(z)=cz/(z+s)$ for some $c,s>0$.
\end{enumerate}
If the third condition fails, then the family in
Section~\ref{sec:family} satisfies
\[
 \det f(A_r)[\{1,4\},\{1,4\}]<0
\]
for arbitrarily large positive integers $r$. Thus the four-dimensional
counterexamples have rational entries and determinant one.
Counterexamples also exist in
$\TP_n\cap\mathrm{SL}_n(\mathbb R)$ for every $n\ge4$.
\end{theorem}

\begin{corollary}\label{cor:TPunbounded}
Let $n\ge4$, and let $f$ be a complete Bernstein function which is
unbounded on the positive real axis. The following are equivalent:
\begin{enumerate}
\item $f(A)\in\mathcal P_n$ for every $A\in\mathcal P_n$;
\item $f(A)\in\mathcal P_n$ for every $A\in\TP_n$;
\item $f(z)=a+bz$, where $a\ge0$ and $b>0$.
\end{enumerate}
\end{corollary}

\begin{corollary}\label{cor:universal}
The nonzero complete Bernstein functions which preserve real
P-matrices in every finite order are exactly
\[
 f(z)=a+bz,\qquad a,b\ge0,\quad a+b>0,
\]
and
\[
 f(z)=a+c\frac{z}{z+s},\qquad a\ge0,\quad c,s>0.
\]
\end{corollary}

\begin{theorem}\label{thm:constant}
For every positive integer $n$ and every $0<\ell<u$, there is a
constant $c=c(n,u/\ell)>0$ with the following property.
Let $\mu$ be a finite positive measure supported on $[\ell,u]$ and put
\[
 g(z)=\int_{[\ell,u]}\frac{z}{z+s}\,\mu(ds),\qquad
 C=\mu([\ell,u]).
\]
If $a>cC$, then $a+g$ preserves every real P-matrix of order $n$.
\end{theorem}

For a positive measure $\mu$ on $(0,\infty)$, write
\[
 M_t=\int_{(0,\infty)}s^t\,\mu(ds),\qquad t\in\mathbb R,
\]
allowing the value $\infty$. For $n\ge4$, put
\begin{equation}\label{eq:rn}
 r_n=\lfloor n/2\rfloor-1.
\end{equation}

\begin{theorem}\label{thm:repair}
Fix $n\ge4$ and let $\mu$ be a finite nonnegative measure on $(0,\infty)$.
Put $g(z)=\int z/(z+s)\,\mu(ds)$.
There is a finite $a\ge0$ such that $aI+g(A)\in\Pn$ for every
$A\in\Pn$ if and only if $M_{r_n},M_{-r_n}<\infty$.
There is a constant $c_n>0$, depending only on $n$, such that every
\[
 a>c_n\sqrt{M_{r_n}M_{-r_n}}
\]
is sufficient. Conversely, any $a\ge0$ with the stated preservation
property satisfies
\[
 \sqrt{M_{r_n}M_{-r_n}}\le2a+M_0.
\]
\end{theorem}

The sufficient bound follows from a simultaneous diagonal scaling
of the whole family of matrix kernels.

\begin{theorem}\label{thm:kernel}
For each $n\ge4$ there is a constant $B_n$ such that every $A\in\Pn$
admits a positive diagonal matrix $Q_A$ satisfying
\[
 \max_{i,j}\left|
 (Q_A^{-1}A(A+sI)^{-1}Q_A)_{ij}\right|
 \le B_n\max(s,s^{-1})^{r_n}\qquad(s>0).
\]
The matrix $Q_A$ is independent of $s$.
\end{theorem}

\begin{proposition}\label{prop:outerlayers}
Let $n\ge4$ and let
\[
 g(z)=bz+\int_{(0,\infty)}\frac{z}{z+s}\,\mu(ds)
\]
be a nonzero complete Bernstein function with $g(0+)=0$.
The following conditions are equivalent:
\begin{enumerate}
\item There is a finite $a\ge0$ such that all principal minors of order
      $2$ of $aI+g(A)$ are positive for every $A\in\Pn$.
\item Either $\mu=0$, or $b=0$ and $M_{r_n},M_{-r_n}<\infty$.
\item There is a finite $a\ge0$ such that $aI+g(A)\in\Pn$ for every
      $A\in\Pn$.
\end{enumerate}
\end{proposition}

\begin{corollary}\label{cor:dimensions}
Let $g(z)=bz+\int_{(0,\infty)}z/(z+s)\,\mu(ds)$, where
$b$ and $\mu$ satisfy the conditions in \eqref{eq:kernel}.
There are finite constants $a_n\ge0$, one for each positive integer
$n$, such that $a_n+g$ preserves $\Pn$ for all $n$, if and only if
either $\mu=0$, or $b=0$ and
\[
 M_j,M_{-j}<\infty\qquad\text{for every integer }j\ge1.
\]
\end{corollary}

\section{An integer matrix with totally positive powers}\label{sec:integer}

A matrix is totally nonnegative if all its nonempty minors are
nonnegative. For a matrix $Y$ of order $n$, put
\begin{equation}\label{eq:initial}
\Delta_{ij}(Y)=\det Y[\{i-k+1,\ldots,i\},\{j-k+1,\ldots,j\}],
\qquad k=\min(i,j).
\end{equation}
These are the initial minors of $Y$.

\begin{lemma}[{\cite[Theorems~9 and 11]{FZ}}]\label{lem:TP}
A real square matrix is strictly totally positive if and only if all its
initial minors are positive. Every invertible totally nonnegative matrix
is the limit of strictly totally positive matrices.
\end{lemma}

\begin{proposition}\label{prop:integer}
There is an entrywise positive integer matrix $X$ of order four with
distinct positive eigenvalues such that
\[
X^k\in\TP_4\quad\text{for every integer }k\ge2,
\qquad \det X[\{1,4\},\{1,4\}]=-26.
\]
Thus, for every integer $m\ge2$, a strictly totally positive integer matrix
has an entrywise positive integer principal $m$th root which is not a
P-matrix.
\end{proposition}

\begin{proof}
Take
\begin{equation}\label{eq:X}
X=\begin{pmatrix}
10&47&14&6\\41&203&113&14\\74&367&223&34\\26&129&80&13
\end{pmatrix}.
\end{equation}
Its entries are positive and its indicated principal minor is
$10\cdot13-6\cdot26=-26$. Direct multiplication and evaluation of
\eqref{eq:initial} give
\[
(\Delta_{ij}(X^2))=
\begin{pmatrix}
3219&15923&9053&1272\\
17459&163790&11178047&998632\\
33173&294105&5018885&147463434\\
11807&48916&60929&243049
\end{pmatrix}.
\]
The corresponding table for $X^3$ is
\[
(\Delta_{ij}(X^3))=
\begin{pmatrix}
1388027&6870201&3964944&566574\\
7590003&835292905&62894938215&7076623086\\
14444676&1720059021&18527950505&588501795498\\
5143236&289126728&292027767&119823157
\end{pmatrix}.
\]
All entries in both tables are positive, so Lemma~\ref{lem:TP} gives
$X^2,X^3\in\TP_4$.
The Cauchy--Binet formula shows that a product of strictly totally positive
matrices is strictly totally positive. Every integer $k\ge2$ is either
$2j$ with $j\ge1$ or $3+2j$ with $j\ge0$, which proves the assertion
about integer powers.

The characteristic polynomial of $X$ is
\[
p_X(x)=x^4-449x^3+6081x^2-4640x+493.
\]
Its values at $0,1/4,3,15,575$, respectively, are
\[
493,\qquad-\frac{75251}{256},\qquad29260,\qquad-165632,
\qquad25961644368.
\]
The intermediate value theorem gives one root in each of the four
successive positive intervals. These exhaust the roots, so the spectrum
is positive and simple. Diagonalization therefore gives
$(X^m)^{1/m}=X$ on the principal branch. Taking $A=X^m$ proves the final
assertion.
\end{proof}

Kushel~\cite[Theorem~7]{Kushel} characterized matrices whose sufficiently
large integer powers are strictly totally positive. The matrix above
has this property from the second power onward, while its negative
principal minor gives the stated counterexamples for every root order.

For $A=X^2$ in Proposition~\ref{prop:integer}, every positive integer
power of $A$ is totally positive, whereas $A^{1/2}=X$ is not a P-matrix.
This gives a counterexample to the real-power statement for PM-matrices
in \cite[Theorem~4.1, p.~363]{MST}.

The proof of Proposition~3.13 in \cite[p.~362]{MST} has a gap.
For a fixed index set $\alpha$ satisfying
$A[\alpha^c,\alpha]=0$, the proof obtains
\[
 f(A)[\alpha,\alpha]=f\bigl(A/\alpha^c\bigr),
\]
where $A/\alpha^c$ denotes the Schur complement of
$A[\alpha^c,\alpha^c]$. It then concludes positivity of all principal
minors, although the zero-block hypothesis was imposed only for the
fixed set $\alpha$. The embedding $\operatorname{diag}(A,1)$ used in
Theorem~3.14 (p.~362) does not supply the missing identities for
the proper principal submatrices of $f(A)$.

\section{A totally positive family of order four}\label{sec:family}

For $r>0$, define
\[
L_r=\begin{pmatrix}
1&0&0&0\\r&1&0&0\\r^2&r&1&0\\r^3&r^2&r&1
\end{pmatrix},\qquad D_r=\diag(r,1,1,r^{-1}),
\]
\[
U_r=\begin{pmatrix}
1&r+1+r^{-1}&r+1+r^{-2}&r^{-2}\\
0&1&r+r^{-1}&r^{-1}\\0&0&1&1\\0&0&0&1
\end{pmatrix},\qquad
E_r=\begin{pmatrix}
1&0&0&0\\2r^{-3}&1&0&0\\r^{-6}&r^{-3}&1&0\\0&0&0&1
\end{pmatrix},
\]
and put
\begin{equation}\label{eq:Ar}
A_r=E_rL_rD_rU_r.
\end{equation}
The determinants of the factors give $\det A_r=1$. Its initial minors are
\begin{equation}\label{eq:table}
(\Delta_{ij}(A_r))=
\begin{pmatrix}
r&r^2+r+1&r^2+r+r^{-1}&r^{-1}\\
r^2+2r^{-2}&r&r^3+r+1&r\\
r^3+r^{-1}+r^{-5}&2r^{-1}+r^{-5}&r&r^3\\
r^4&r^{-3}&r^{-4}&1
\end{pmatrix}.
\end{equation}
They are positive for every $r>0$, so Lemma~\ref{lem:TP} proves
$A_r\in\TP_4$. The first and fourth rows of $E_r$ coincide with those of
the identity; in particular,
\begin{equation}\label{eq:f1}
\det A_r[\{1,4\},\{1,4\}]=2r^2+1.
\end{equation}

Write
\begin{equation}\label{eq:char}
p_r(x)=\det(xI-A_r)=x^4-a_rx^3+b_rx^2-c_rx+1.
\end{equation}
Expansion of this determinant gives
\begin{equation}\label{eq:abc}
a_r=r^4+O(r^3),\qquad b_r=r^4+O(r^3),\qquad
c_r=3r^2+O(r).
\end{equation}
The principal minors have the following leading terms, listed in
lexicographic order of their index sets:
\[
\begin{array}{c|l}
\text{order}&\text{leading terms}\\\hline
1&r,\ r^3,\ r^4,\ r^2\\
2&r,\ r^3,\ 2r^2,\ r^3,\ r^4,\ r^3\\
3&r,\ r^2,\ r^2,\ r^2
\end{array}.
\]
Each remainder is $O(r^{-1})$ times the displayed term.
Summing these minors gives \eqref{eq:abc}.

\begin{lemma}\label{lem:spec}
Put $\alpha=(3+\sqrt5)/2$ and $\beta=(3-\sqrt5)/2$.
For sufficiently large $r$, the eigenvalues of $A_r$, in decreasing order,
satisfy
\begin{equation}\label{eq:spec}
\begin{aligned}
\lambda_1&=r^4(1+O(r^{-1})),&
\lambda_2&=1+O(r^{-1}),\\
\lambda_3&=\alpha r^{-2}(1+O(r^{-1})),&
\lambda_4&=\beta r^{-2}(1+O(r^{-1})).
\end{aligned}
\end{equation}
They are distinct and positive.
\end{lemma}

\begin{proof}
On every bounded real $z$-interval, \eqref{eq:abc} gives the convergences
\[
r^{-16}p_r(r^4z)\longrightarrow z^3(z-1),\qquad
r^{-4}p_r(z)\longrightarrow z^2(1-z),\qquad
p_r(r^{-2}z)\longrightarrow z^2-3z+1.
\]
The coefficient errors are $O(r^{-1})$ in each case. Choose short intervals around the simple root $1$ of each of the first
two limiting polynomials and around the simple roots $\alpha,\beta$
of the third, with endpoints on opposite sides of each root.
For sufficiently large $r$, sign changes give a root in each interval.
The four scaled intervals are disjoint and contain all roots of $p_r$.
On these intervals the limiting derivatives stay bounded away from zero;
the same holds for the derivatives of the normalized polynomials.
The mean value theorem then bounds each root displacement by $O(r^{-1})$,
which proves \eqref{eq:spec}.
\end{proof}

Write $ij$ for the increasing index set $\{i,j\}$.
For $I\in\{14,23\}$, put
\[
 F_I(t)=\det A_r^t[I,I].
\]
For sufficiently large $r$, the eigenvalues are positive and distinct.
Choose a spectral decomposition
$A_r=V\operatorname{diag}(\lambda_1,\ldots,\lambda_4)V^{-1}$ and put
\[
 w_{ij}^I=\det V[I,ij]\det V^{-1}[ij,I].
\]
The Cauchy--Binet formula gives
\begin{equation}\label{eq:sum}
 F_I(t)=\sum_{i<j}w_{ij}^I(\lambda_i\lambda_j)^t.
\end{equation}

\begin{lemma}\label{lem:filter}
Let $\{i,j,k,l\}=\{1,2,3,4\}$, with $i<j$ and $k<l$. Then
\begin{equation}\label{filterformula}
 w_{ij}^I=
 \frac{\det\!\left(((A_r-\lambda_k I_4)(A_r-\lambda_l I_4))[I,I]\right)}
 {(\lambda_i-\lambda_k)(\lambda_i-\lambda_l)
  (\lambda_j-\lambda_k)(\lambda_j-\lambda_l)}.
\end{equation}
\end{lemma}

\begin{proof}
Set $q(x)=(x-\lambda_k)(x-\lambda_l)$.
Apply Cauchy--Binet to
$q(A_r)=V\operatorname{diag}(q(\lambda_1),\ldots,q(\lambda_4))V^{-1}$.
In its principal minor on $I$, every summand vanishes except the one
indexed by $ij$, which is $w_{ij}^I q(\lambda_i)q(\lambda_j)$.
This is \eqref{filterformula}.
\end{proof}

To evaluate the numerator in \eqref{filterformula}, set
$M_I=A_r[I,I]$ and $N_I=A_r^2[I,I]$.
For scalars $s,p$,
\begin{equation}\label{quadratic}
 \det(N_I-sM_I+pI_2)
 =\det N_I-sH_I+s^2\det M_I
   +p\operatorname{tr}N_I-ps\operatorname{tr}M_I+p^2,
\end{equation}
where
\begin{equation}\label{mixed}
 H_I=\operatorname{tr}(A_r)\det A_r[I,I]
       -\sum_{k\notin I}\det A_r[I\cup\{k\},I\cup\{k\}].
\end{equation}
Indeed, expansion of the two by two determinant identifies $H_I$ as
$(M_I)_{11}(N_I)_{22}+(M_I)_{22}(N_I)_{11}
 -(M_I)_{12}(N_I)_{21}-(M_I)_{21}(N_I)_{12}$.
Expanding the two three by three determinants in \eqref{mixed} gives
exactly this expression.

The principal-minor estimates above and Cauchy--Binet give
\begin{equation}\label{smalltable}
\begin{array}{c|ccc}
 I&\det M_I&H_I&\det N_I\\\hline
 14&2r^2+1&O(r^6)&r^7+O(r^6)\\
 23&r^3+O(r^2)&r^7+O(r^6)&2r^7+O(r^6)
\end{array}.
\end{equation}
Here the estimates for $H_I$ follow from \eqref{mixed}, because the
three by three principal minors are positive and sum to
$c_r=3r^2+O(r)$. For the last column, write
\[
 \det N_I=\sum_{|J|=2}\det A_r[I,J]\det A_r[J,I].
\]
Direct multiplication of the factors defining $A_r$ gives
\begin{align*}
 \det A_r[14,23]\det A_r[23,14]&=r^7+O(r^6),\\
 \det A_r[23,24]\det A_r[24,23]&=r^7+O(r^6).
\end{align*}
All other summands for $I=14$ are $O(r^6)$; for $I=23$, only the
summands $J=14,24$ have order $r^7$.
Finally, \eqref{eq:spec} gives
$\operatorname{tr}A_r=O(r^4)$ and
$\operatorname{tr}(A_r^2)=O(r^8)$.
Since both matrices have positive diagonal entries, it follows that
$\operatorname{tr}M_I=O(r^4)$ and $\operatorname{tr}N_I=O(r^8)$.
Consequently, uniformly for $s=1+O(r^{-1})$ and $p=O(r^{-2})$,
\begin{equation}\label{filterestimate}
 \det(N_I-sM_I+pI_2)=r^7(1+O(r^{-1}))
 \qquad(I=14,23).
\end{equation}
For $I=23$, the two terms of order $r^7$ in \eqref{quadratic} have
coefficients $2$ and $-1$; all other terms are $O(r^6)$.

\begin{lemma}\label{lem:w}
As $r\to\infty$, the coefficients $w_{ij}=w_{ij}^{14}$ in
\eqref{eq:sum} satisfy
\begin{equation}\label{eq:w}
\begin{array}{c|c}
 ij&w_{ij}\\\hline
12&r^{-1}(1+O(r^{-1}))\\
13&-r/\sqrt5\,(1+O(r^{-1}))\\
14&r/\sqrt5\,(1+O(r^{-1}))\\
23&r/\sqrt5\,(1+O(r^{-1}))\\
24&-r/\sqrt5\,(1+O(r^{-1}))\\
34&2r^{-1}(1+O(r^{-1}))
\end{array}.
\end{equation}
In addition,
\begin{equation}\label{eq:wsums}
 w_{13}+w_{14}=-r^{-1}+O(r^{-2}),\qquad
 w_{23}+w_{24}=1-2r^{-1}+O(r^{-2}).
\end{equation}
\end{lemma}

\begin{proof}
For $ij=13$ or $14$, the complementary indices in
\eqref{filterformula} contain $2$ and one of $3,4$.
Thus $s=\lambda_k+\lambda_l=1+O(r^{-1})$ and
$p=\lambda_k\lambda_l=O(r^{-2})$.
The spectral estimates give, respectively,
\[
 (\lambda_1-\lambda_2)(\lambda_1-\lambda_4)
 (\lambda_3-\lambda_2)(\lambda_3-\lambda_4)
 =-\sqrt5\,r^6(1+O(r^{-1}))
\]
and
\[
 (\lambda_1-\lambda_2)(\lambda_1-\lambda_3)
 (\lambda_4-\lambda_2)(\lambda_4-\lambda_3)
 =\sqrt5\,r^6(1+O(r^{-1})).
\]
Equations \eqref{filterformula} and \eqref{filterestimate} prove
$w_{13}^I=-r/\sqrt5\,(1+O(r^{-1}))$ and
$w_{14}^I=r/\sqrt5\,(1+O(r^{-1}))$ for both $I=14,23$.

Since $\det A_r^t=1$, \eqref{eq:jprincipal} gives
$F_{14}(-t)=F_{23}(t)$.
Also $(\lambda_i\lambda_j)^{-1}=\lambda_k\lambda_l$ for complementary
pairs. The six products are distinct, so comparison in \eqref{eq:sum}
gives $w_{23}^{14}=w_{14}^{23}$ and
$w_{24}^{14}=w_{13}^{23}$.
This proves the four middle entries of the table.

For the two extreme weights, put
$R=\lambda_1\lambda_2=r^4(1+O(r^{-1}))$, so
$\lambda_3\lambda_4=R^{-1}$.
The four middle terms in \eqref{eq:sum} have sum
$-r^3(1+O(r^{-1}))$ at both $t=1$ and $t=-1$.
The principal minors of $A_r$ give
\[
 F_{14}(1)=2r^2+1,\qquad
 F_{14}(-1)=F_{23}(1)=r^3+O(r^2).
\]
The remaining two entries follow by inverting the system
\[
 \begin{pmatrix}R&R^{-1}\\R^{-1}&R\end{pmatrix}
 \begin{pmatrix}w_{12}\\w_{34}\end{pmatrix}
 =\begin{pmatrix}r^3(1+O(r^{-1}))\\2r^3(1+O(r^{-1}))\end{pmatrix}.
\]

To prove \eqref{eq:wsums}, let $P_1$ be the spectral projector associated
with $\lambda_1$. Taking the residue of the diagonal resolvent gives
\[
 (P_1)_{ii}
 =\frac{\det(\lambda_1I_3-A_r[\widehat i,\widehat i])}
        {p_r'(\lambda_1)},\qquad
 \widehat i=[4]\setminus\{i\}.
\]
For $i=1,4$, the remaining coefficients of the numerator's cubic
polynomial are bounded by $b_r=O(r^4)$ and $c_r=O(r^2)$, since they
are sums of principal minors of $A_r$. Moreover,
$\lambda_1-a_r=-(\lambda_2+\lambda_3+\lambda_4)=O(1)$ and
$(A_r)_{ii}=O(r^2)$. Consequently
\[
 \det(\lambda_1I_3-A_r[\widehat i,\widehat i])
 =\lambda_1^2(\lambda_1-a_r+(A_r)_{ii})
   +O(r^4\lambda_1+r^2)=O(r^{10}).
\]
As $p_r'(\lambda_1)=r^{12}(1+O(r^{-1}))$, we obtain
$(P_1)_{11}+(P_1)_{44}=O(r^{-2})$.
Since $P_1$ has rank one, comparison of the coefficient of $x$ in
$\det((I_4+xP_1)[14,14])$, using Cauchy--Binet, gives
\[
 w_{12}+w_{13}+w_{14}=(P_1)_{11}+(P_1)_{44}.
\]
Together with $\sum_{i<j}w_{ij}=1$ and the extreme weights in
\eqref{eq:w}, this proves \eqref{eq:wsums}.
\end{proof}

Put $f_r(t)=F_{14}(t)$ and $\rho_{ij}=\lambda_i\lambda_j$.
Then \eqref{eq:sum} is the six-term expansion
$f_r(t)=\sum_{i<j}w_{ij}\rho_{ij}^t$.

\section{Negative minors on compact exponent intervals}\label{sec:negative}

\begin{proposition}\label{prop:limit}
For $0<|t|<1$, the family \eqref{eq:Ar} satisfies
\begin{equation}\label{eq:limit}
\frac{f_r(t)}{r^{1+2|t|}}
\longrightarrow-\frac{\alpha^{|t|}-\beta^{|t|}}{\sqrt5}<0.
\end{equation}
The convergence is uniform on every compact subset of
$(-1,0)\cup(0,1)$.
\end{proposition}

\begin{proof}
For fixed $0<t<1$, the six terms of \eqref{eq:sum}, in the order
$12,13,14,23,24,34$, have respective orders
\[
r^{4t-1},\quad r^{1+2t},\quad r^{1+2t},\quad
r^{1-2t},\quad r^{1-2t},\quad r^{-1-4t}.
\]
Since $1+2t>4t-1$ and $1+2t>1-2t$, the terms with indices $13$ and
$14$ give \eqref{eq:limit}. For $t=-s$, $0<s<1$, the dominant terms
are those with indices $23$ and $24$; their sum has the same negative
limit because $\alpha\beta=1$.

More precisely, for $K=[-b,-a]\cup[a,b]$, $0<a\le b<1$, the relative
errors in \eqref{eq:spec} and \eqref{eq:w} remain $O_K(r^{-1})$ after
raising the spectral factors to the bounded powers in $K$. The extreme
spectral term is smaller than the dominant terms by at least
$r^{-2(1-b)}$, and the opposite intermediate pair is smaller by at least
$r^{-4a}$. Hence
\begin{equation}\label{eq:uniform}
\frac{f_r(t)}{r^{1+2|t|}}
=-\frac{\alpha^{|t|}-\beta^{|t|}}{\sqrt5}
+O_K\bigl(r^{-1}+r^{-2(1-b)}+r^{-4a}\bigr).
\end{equation}
\end{proof}

\begin{proof}[Proof of Theorem~\ref{thm:uniform}]
The limiting function in \eqref{eq:limit} is bounded above by a strictly
negative constant on $K$. Thus a single sufficiently large $r$ gives
\eqref{eq:negative} for the whole of $K$. The matrices $A_r$ are strictly
totally positive and have determinant one. Taking $r$ to be an integer also makes all entries rational.

For $n>4$, take such an $A=A_r$ and put $D=\diag(A,I_{n-4})$.
This is invertible and totally nonnegative, and its $14$ principal minor
at each power equals that of $A^t$. By Lemma~\ref{lem:TP}, $D$ is a limit
of matrices in $\TP_n$.
For $Y$ sufficiently close to $D$, its spectrum still has no nonpositive
real eigenvalues. The integral formulas~\eqref{eq:int1}--\eqref{eq:int2}
show that $(Y,t)\mapsto Y^t$ is jointly continuous there for $t\in K$.
Indeed, the resolvents are uniformly bounded on bounded $s$-intervals
and are uniformly $O(s^{-1})$ at infinity. For $a\le t\le b$, the
integrands in the two formulas are bounded by constant multiples of
$s^{a-1}$ and $s^{-b}$ on $(0,1)$, and of $s^{b-2}$ and $s^{-a-1}$
on $(1,\infty)$, respectively. These bounds are integrable, so dominated
convergence gives the required joint continuity, including at the
repeated eigenvalue $1$ of $D$.
Compactness supplies a uniform strict negative bound for the indicated
minor, so a sufficiently close totally positive matrix retains the bound
for every $t\in K$. Multiplying it by the positive scalar
$(\det Y)^{-1/n}$ normalizes its determinant without changing the sign
of any powered minor.

The same continuity argument shows that the inequalities throughout $K$
persist in a neighborhood of the chosen matrix. Total positivity is an
open condition. Intersection with determinant one proves the relative
openness assertion.
\end{proof}

\begin{proof}[Proof of Theorem~\ref{thm:main}]
Proposition~\ref{prop:low} gives every exponent in $[-1,1]$ for $n\le3$.
The rotation example following its proof excludes every $|t|>1$ for
$n\ge2$. For $n\ge4$, Theorem~\ref{thm:uniform}, with $a=b=|t|$,
excludes every $0<|t|<1$, even for totally positive input.
The exponents $-1,0,1$ preserve every P-matrix in every order, by
inversion, the identity matrix, and the definition. Positive scalars
retain positivity at all real powers, which handles $n=1$.
\end{proof}

Put
\[
 L=\log r,\qquad \ell=\log\alpha,\qquad
 c=\left(\frac{\sqrt5}{8\ell}\right)^{1/2},\qquad
 d=\frac{3\ell}{\sqrt5}.
\]

\begin{proposition}\label{prop:zeros}
For sufficiently large $r$, $f_r$ has exactly four real zeros, all simple,
with
\[
 -1<\xi_1(r)<\xi_2(r)<0<\xi_3(r)<\xi_4(r)<1.
\]
It is negative precisely on
$(\xi_1,\xi_2)\cup(\xi_3,\xi_4)$. As $r\to\infty$,
\begin{equation}\label{eq:zerospeeds}
\begin{aligned}
 \xi_1&=-1+\frac{\log2}{2L}
              +\frac{d\log2}{4L^2}+O(L^{-3}),\\
 \xi_2&=-\frac{c}{\sqrt{rL}}-\frac{c^2}{r}
              +O((rL)^{-1}),\\
 \xi_3&=\frac{c}{\sqrt{rL}}-\frac{c^2}{r}
              +O((rL)^{-1}),\\
 \xi_4&=1-\frac1{rL}
              \left(1+\frac d{2L}+O(L^{-2})\right).
\end{aligned}
\end{equation}
The corresponding local limits, uniformly for $x$ in each bounded interval,
are
\begin{equation}\label{eq:zeroprofiles}
\begin{aligned}
 r^{-3}f_r(-1+x/L)&\longrightarrow 2e^{-4x}-e^{-2x},\\
 f_r\bigl(x/\sqrt{rL}\bigr)&\longrightarrow 1-x^2/c^2,\\
 r^{-2}f_r\bigl(1-x/(rL)\bigr)&\longrightarrow 2-2x.
\end{aligned}
\end{equation}
\end{proposition}

\begin{proof}
A nonzero sum of $m$ exponentials with distinct real exponents has at
most $m-1$ real zeros counted with multiplicity. This follows by factoring
out one exponential, differentiating, and applying Rolle's theorem
inductively. In \eqref{eq:sum}, $m=6$ for sufficiently large $r$.
The extreme weights are positive by \eqref{eq:w}, so $f_r$ is positive
at both ends of the real axis. Its total number of zeros counted with
multiplicity is therefore even and at most four.
We have $f_r(0)=1$, $f_r(1)=2r^2+1$, and
$f_r(-1)=F_{23}(1)=\det A_r[23,23]>0$.
Proposition~\ref{prop:limit} makes $f_r(1/2)$ and $f_r(-1/2)$ negative
for large $r$. The four intervening sign changes give all the zeros,
each simple, and the asserted negative intervals.

Write $u_{ij}=\log(\lambda_i\lambda_j)$. Equations
\eqref{eq:spec} and \eqref{eq:w} give
\[
\begin{array}{ll}
u_{13}=2L+\ell+O(r^{-1}),&u_{14}=2L-\ell+O(r^{-1}),\\
u_{23}=-2L+\ell+O(r^{-1}),&u_{24}=-2L-\ell+O(r^{-1}),\\
u_{12}=4L+O(r^{-1}),&u_{34}=-4L+O(r^{-1}).
\end{array}
\]
Together with \eqref{eq:wsums}, these imply
\begin{align*}
 f'_r(0)
 ={}&2L\bigl[(w_{13}+w_{14})-(w_{23}+w_{24})\bigr]\\
 &+\ell(w_{13}-w_{14}+w_{23}-w_{24})+O(1)
 =-2L+O(1).
\end{align*}
Replacing the weights and exponents in the four middle terms by their
leading terms gives
\[
 S_r(t)=-\frac{4r}{\sqrt5}\sinh(\ell t)\sinh(2Lt).
\]
Set $\eta=(rL)^{-1/2}$ and $\varepsilon=\sqrt{L/r}$.
Then $S_r(\eta x)=-x^2/c^2+O_M(\varepsilon^2)$ uniformly for $|x|\le M$.
For $|t|\le M\eta$, termwise differentiation gives
$(f_r-S_r)''(t)=O_M(L^2)$: the weight errors are $O(1)$,
the exponent errors are $O(r^{-1})$, and $L\eta\to0$.
Since $S_r(0)=S'_r(0)=0$, Taylor's formula and the preceding derivative
estimate give, uniformly for $|x|\le M$,
\begin{equation}\label{eq:centralprofile}
 f_r(\eta x)=1-\frac{x^2}{c^2}-2\varepsilon x
                    +O_M(\eta+\varepsilon^2).
\end{equation}
The two central zeros therefore satisfy
$\xi_2/\eta=-c+O(\varepsilon)$ and
$\xi_3/\eta=c+O(\varepsilon)$.
Substitution back into \eqref{eq:centralprofile} yields
\[
 \frac{\xi_2}{\eta}=-c-c^2\varepsilon
                 +O(\eta+\varepsilon^2),\qquad
 \frac{\xi_3}{\eta}=c-c^2\varepsilon
                 +O(\eta+\varepsilon^2).
\]
Now $\eta\varepsilon=r^{-1}$ and
$\eta(\eta+\varepsilon^2)=O((rL)^{-1})$,
which proves the two central expansions and their local limit.

For the left endpoint put
$C(s)=(\alpha^s-\beta^s)/\sqrt5$.
The three dominant terms of \eqref{eq:sum} give
\begin{equation}\label{eq:leftprofile}
 r^{-3}f_r(-1+x/L)
 =2e^{-4x}-C(1-x/L)e^{-2x}+O_M(r^{-1})
 \qquad(|x|\le M).
\end{equation}
Since $C(1)=1$, this proves the first limit in
\eqref{eq:zeroprofiles}. Its simple zero is $x=(\log2)/2$;
the sign changes on its two sides show that
$\delta:=1+\xi_1=O(L^{-1})$.
At this zero, \eqref{eq:leftprofile} gives
\[
 \log2-2L\delta=\log C(1-\delta)+O(r^{-1})
                     =-d\delta+O(\delta^2)+O(r^{-1}),
\]
because $C'(1)=\ell(\alpha+\beta)/\sqrt5=d$.
Thus $\delta=(\log2)/(2L-d)+O(L^{-3})$, proving the first
expansion in \eqref{eq:zerospeeds}.

At the right endpoint, direct differentiation of \eqref{eq:sum} at $t=1$ gives
\begin{equation}\label{eq:rightderivative}
 f'_r(1)=r^3\left(2L-d+O(L/r)\right).
\end{equation}
Indeed, term $12$ contributes $4Lr^3+O(r^2L)$, and terms $13,14$
together contribute $-r^3(2L+d)+O(r^2L)$; the remaining terms are
$O(L/r)$. The same expansion gives
$f''_r(t)=O_M(r^3L^2)$ for $|t-1|\le M/(rL)$.
Since $f_r(1)=2r^2+1$, Taylor's formula proves the last limit in
\eqref{eq:zeroprofiles} and locates $1-\xi_4=O((rL)^{-1})$.
Applying it again at the zero and using \eqref{eq:rightderivative}
gives
\[
 1-\xi_4=\frac{f_r(1)}{f'_r(1)}\bigl(1+O(r^{-1})\bigr)
 =\frac1{rL}\left(1+\frac d{2L}+O(L^{-2})\right).
\]
\end{proof}

\section{Functions represented by positive kernels}\label{sec:functions}

\begin{proof}[Proof of Theorem~\ref{thm:functions}]
Positive constant functions preserve P-matrices, so assume that $f$ is
nonconstant.
For nonreal $z$,
\begin{equation}\label{eq:kernelimag}
 \frac{\operatorname{Im}f(z)}{\operatorname{Im}z}
 =b+\int_{(0,\infty)}\frac{s}{|z+s|^2}\,\mu(ds)>0.
\end{equation}
Thus $f$ has no nonreal zero, and $f(x)>0$ for $x>0$.
The spectrum of a P-matrix avoids $(-\infty,0]$ by \eqref{eq:p},
so $f(A)$ is invertible. Each positive real eigenvalue $\lambda$ of $A$
contributes $f(\lambda)>0$ to its determinant, and each nonreal conjugate
pair contributes $|f(\lambda)|^2>0$. Hence $\det f(A)>0$.
Equation~\eqref{eq:res2} and the matrix integral
\eqref{eq:matrixkernel} give positive diagonal entries.

For the complementary principal minors, first consider a nonzero finite sum
\[
 g(z)=a+bz+\sum_{j=1}^m c_j\frac{z}{z+s_j},\qquad
 a,b\ge0,\quad c_j,s_j>0,
\]
combining terms with equal $s_j$. The sum may be empty if $a+b>0$.
If $g$ is nonconstant, its zeros are real by \eqref{eq:kernelimag}, and
\[
 g'(x)=b+\sum_{j=1}^m\frac{c_js_j}{(x+s_j)^2}>0
\]
away from its poles. Every zero is therefore simple and nonpositive.
Partial fractions give
\begin{equation}\label{eq:reciprocalfinite}
 \frac1{g(z)}=a_0+\frac{d_0}{z}
       +\sum_{\ell}\frac{d_\ell}{z+u_\ell},\qquad
 a_0,d_0\ge0,\quad d_\ell,u_\ell>0.
\end{equation}
Indeed, the residue at a zero $-u_\ell$ is
$1/g'(-u_\ell)>0$. If $g(0)=0$, its residue at zero is
$d_0=1/g'(0)>0$; otherwise $d_0=0$.
The constant term is $a_0=\lim_{x\to\infty}1/g(x)\ge0$.
For a positive constant $g$, the expansion has only the constant term.

For each fixed $A\in\Pn$, truncate the measure in \eqref{eq:kernel}
to compact subintervals of $(0,\infty)$ and approximate it there by
finite positive atomic measures. The bounded kernel near zero and its
$O(s^{-1})$ decay at infinity show that the resulting nonzero finite
sums $f_k$ can be chosen with
\[
 f_k(1)\longrightarrow f(1),\qquad f_k(A)\longrightarrow f(A).
\]
Since $f(A)$ is invertible, their inverses converge as well.
For $i\in[n]$, put
\[
 m_i=\inf_{s\ge0}(1+s)((A+sI)^{-1})_{ii}>0.
\]
The expression is continuous and positive for finite $s$, by the
diagonal cofactor formula, and tends to $1$ at infinity. Thus
$0<m_i\le1$. Applying \eqref{eq:reciprocalfinite} to $f_k$ gives
\[
 (f_k(A)^{-1})_{ii}\ge\frac{m_i}{f_k(1)}.
\]
Passing to the limit yields
\begin{equation}\label{eq:inversediagonal}
 (f(A)^{-1})_{ii}\ge\frac{m_i}{f(1)}>0.
\end{equation}
Taking $B=f(A)$ in \eqref{eq:jprincipal} gives
\[
 \det f(A)[\widehat i,\widehat i]
 =\det f(A)\,(f(A)^{-1})_{ii}>0,
 \qquad \widehat i=[n]\setminus\{i\}.
\]
For $n\le3$, these orders exhaust all nonempty principal minors.

Positive scalar multiplication, addition of a nonnegative multiple of
the identity, and inversion preserve P-matrices. Since
\[
 A(A+sI)^{-1}=(I+sA^{-1})^{-1},
\]
the two functions in the third condition preserve P-matrices in every
order. Now assume $f(0+)=0$. For the family $A_r$, the spectral
decomposition and the Cauchy--Binet formula give
\begin{equation}\label{eq:functionminor}
 \det f(A_r)[14,14]
 =\sum_{i<j}w_{ij}f(\lambda_i)f(\lambda_j),
\end{equation}
with the eigenvalues and weights in \eqref{eq:spec} and \eqref{eq:w}.

First suppose that
\[
 D:=f'(0+)=b+\int_{(0,\infty)}\frac{\mu(ds)}s<\infty,
 \qquad F:=f(1).
\]
Then $f(x)=Dx+o(x)$ at zero. The four middle terms in
\eqref{eq:functionminor}, together with the extreme terms, give
\begin{equation}\label{eq:finiteD}
 r\det f(A_r)[14,14]
 =f(\lambda_1)(F-D)+DF+o(f(\lambda_1)+1).
\end{equation}
Here the terms $12$, $13+14$, and $23+24$ have respective leading
contributions $f(\lambda_1)F$, $-f(\lambda_1)D$, and $DF$ after
multiplication by $r$, and the $34$ term tends to zero.

If $C:=f(+\infty)<\infty$, then $b=0$ and $C=\mu((0,\infty))$.
It follows that
\begin{equation}\label{eq:boundedkernel}
 r\det f(A_r)[14,14]\longrightarrow (C+D)F-CD.
\end{equation}
For the probability measure $\mu/C$, apply strict concavity of
$u/(1+u)$ to $u=1/s$. The means of $u$ and $u/(1+u)$ are
$D/C$ and $F/C$, respectively, so Jensen's inequality gives
\[
 \frac FC\le\frac{D/C}{1+D/C}.
\]
The limit in \eqref{eq:boundedkernel} is nonpositive. Equality holds
only when $\mu$ is supported at one point, giving $f(z)=Cz/(z+s)$;
otherwise the limit is strictly negative.
If $f(+\infty)=\infty$, equation~\eqref{eq:finiteD} gives
\begin{equation}\label{eq:unboundedkernel}
 \frac{r\det f(A_r)[14,14]}{f(\lambda_1)}
 \longrightarrow F-D
 =-\int_{(0,\infty)}\frac{\mu(ds)}{s(1+s)}<0,
\end{equation}
unless $\mu=0$, which is the previously treated linear case.
For finite $D$, therefore, the negative minor occurs for every
sufficiently large $r$ whenever the third condition fails.

Now suppose that $f'(0+)=\infty$. Put
\[
 x=r^{-2},\qquad h=f(\beta x),\qquad
 \Delta=f(\alpha x)-f(\beta x)>0.
\]
The integral representation gives
\begin{equation}\label{eq:incrementderivative}
 \frac\Delta x=(\alpha-\beta)
 \left[b+\int_{(0,\infty)}
 \frac{s}{(s+\alpha x)(s+\beta x)}\,\mu(ds)\right]
 \longrightarrow\infty
\end{equation}
by monotone convergence. We also have
\begin{equation}\label{eq:incrementsubsequence}
 \limsup_{x\downarrow0}\frac{\Delta}{\sqrt{x}\,h}=\infty.
\end{equation}
Otherwise, with $q=\alpha/\beta>1$ and $y=\beta x$, some $C>0$
would satisfy
\[
 f(qy)\le(1+C\sqrt y)f(y)
\]
for all sufficiently small $y$. For a fixed small $y_0>0$, iteration
would imply
\[
 f(y_0q^{-k})\ge
 f(y_0)\prod_{j=1}^k(1+C\sqrt{y_0}\,q^{-j/2})^{-1}.
\]
The product on the right has a positive limit, contradicting $f(0+)=0$.
Choose $x_k\downarrow0$ along which the quotient in
\eqref{eq:incrementsubsequence} tends to infinity, and put
$r_k=x_k^{-1/2}$.

The function $f$ is increasing and concave with $f(0+)=0$, so
\[
 \min(1,c)f(y)\le f(cy)\le\max(1,c)f(y)\qquad(c,y>0).
\]
Along the chosen sequence this gives
\[
 g:=f(\lambda_4)=h(1+O(r^{-1})),\qquad
 \delta:=f(\lambda_3)-f(\lambda_4)
       =\Delta+O(h/r)=\Delta(1+o(1)).
\]
Put $u=f(\lambda_1)$ and $v=f(\lambda_2)$. The exact grouping
\[
\begin{split}
 \det f(A_r)[14,14]
 ={}&(uw_{13}+vw_{23})\delta
    +[u(w_{13}+w_{14})+v(w_{23}+w_{24})]g\\
    &+w_{12}uv+w_{34}f(\lambda_3)g
\end{split}
\]
and \eqref{eq:w}--\eqref{eq:wsums} yield
\begin{equation}\label{eq:infiniteDerror}
 \det f(A_r)[14,14]
 =-\frac r{\sqrt5}(u-v)\delta
   +O(u\delta+uh/r+h+u/r).
\end{equation}
Here $v=O(1)$, $u$ is bounded below, and $h\to0$.
Strict increase of $f$ gives
$\liminf_{r\to\infty}(u-v)/u>0$.
Each error is $o(ru\delta)$: besides $1/r\to0$, this follows from
\[
 \frac{h}{r^2\delta}\longrightarrow0,\qquad
 \frac{h}{r\delta}\longrightarrow0,\qquad
 \frac1{r^2\delta}\longrightarrow0,
\]
using $\delta\sim\Delta$, $\Delta/x\to\infty$, and
$\Delta/(\sqrt{x}\,h)\to\infty$.
Thus the minor is negative along the chosen sequence.
Replacing each $r_k$ by a nearest positive integer changes $x_k$ by a
relative $O(r_k^{-1})$. The same concavity bound changes $h$ by
$O(h/r_k)$ and $\Delta$ by $O(h/r_k)=o(\Delta)$, so the subsequence
condition and the strict negative sign persist. Passing to an increasing
integer subsequence proves the rational determinant-one assertion.

For $n>4$, fix one of these four-dimensional counterexamples and set
$D_0=\diag(A_r,I_{n-4})$. This matrix is nonsingular and
totally nonnegative, and its $f$-image has the same negative $14$ minor.
By Lemma~\ref{lem:TP}, it is a limit of strictly totally positive
matrices. The uniform resolvent bounds used in the proof of
Theorem~\ref{thm:uniform} give, for $Y$ near $D_0$,
\[
 \|Y(Y+sI)^{-1}\|\le\frac{C}{1+s}\qquad(s>0).
\]
The measure condition in \eqref{eq:kernel} and dominated convergence
therefore show that $f(Y)\to f(D_0)$ as $Y\to D_0$, including at its
repeated eigenvalue $1$. A sufficiently close strictly totally positive
approximation retains the negative minor. Since $\det D_0=1$,
normalizing these approximations by the positive scalar
$(\det Y)^{-1/n}$ still gives convergence to $D_0$ and hence preserves
the negative minor for sufficiently close approximations.
\end{proof}

For example, if $f(z)=\sum_{j=1}^m c_jz/(z+s_j)$ with $c_j,s_j>0$,
the limit in \eqref{eq:boundedkernel} becomes
\[
 \lim_{r\to\infty}r\det f(A_r)[14,14]
 =-\sum_{i<j}
 \frac{c_ic_j(s_i-s_j)^2}{s_is_j(1+s_i)(1+s_j)}.
\]
Thus two kernels with different parameters already give a strictly
negative limit. For $f(z)=z/(z+1)+z/(z+2)$, its value is $-1/12$.

\begin{proof}[Proof of Corollary~\ref{cor:TPunbounded}]
Put $f^{\#}(z)=1/f(1/z)$. To obtain its kernel representation,
approximate $f$ by nonzero finite positive kernel sums $f_k$,
retaining its constant and linear terms, so that $f_k(x)\to f(x)$
for $x>0$. Substituting $1/z$ in \eqref{eq:reciprocalfinite} shows
that each $f_k^{\#}$ has the form \eqref{eq:kernel}. Put
\[
 \Phi_u(z)=\frac{z}{(1-u)z+u},\qquad 0\le u\le1.
\]
The endpoint values are $\Phi_0(z)=1$ and $\Phi_1(z)=z$;
for $s=u/(1-u)$, the interior values are
$(1+s)z/(z+s)$. Hence
$f_k^{\#}(z)=\int_{[0,1]}\Phi_u(z)\,\eta_k(du)$ for a finite
positive measure with total mass $f_k^{\#}(1)=1/f_k(1)$.
These masses are bounded, so a subsequence converges weakly to a
positive measure $\eta$. Continuity in $u$ gives
\[
 \frac1{f(1/x)}=\int_{[0,1]}\Phi_u(x)\,\eta(du),\qquad x>0.
\]
The endpoint masses give the constant and linear terms. Push forward
$(1-u)^{-1}\eta(du)$ on $(0,1)$ under $s=u/(1-u)$ to obtain the
kernel measure $\mu^{\#}$. It satisfies
$\int(1+s)^{-1}\mu^{\#}(ds)=\eta((0,1))<\infty$.
Both sides are holomorphic off $(-\infty,0]$, so the identity holds
there as well. Thus $f^{\#}$ is a nonzero complete Bernstein function.
Unboundedness of $f$ implies $f^{\#}(0+)=0$.

Suppose that the second condition holds. Put
$D=\diag(1,-1,\ldots,(-1)^{n-1})$. For $A\in\TP_n$ and index
sets $I,J$ of the same size, the signs in \eqref{eq:jacobi} cancel to give
\[
 \det(DA^{-1}D)[I,J]
 =\frac{\det A[J^c,I^c]}{\det A}>0.
\]
Thus $B=DA^{-1}D$ is strictly totally positive. By the assumed
preservation, $f(B)$ is a P-matrix. Inversion and diagonal sign
similarity preserve all principal minors, and functional calculus gives
\[
 f^{\#}(A)=f(A^{-1})^{-1}=D f(B)^{-1}D\in\mathcal P_n.
\]
The zero-origin classification in Theorem~\ref{thm:functions} now
implies that $f^{\#}(z)=cz$ or $f^{\#}(z)=cz/(z+s)$, with all
displayed parameters positive. Taking the reciprocal dual again gives
$f(z)=z/c$ or $f(z)=(1+sz)/c$, respectively. These are precisely the
affine functions in the third condition.

Conversely, if $a\ge0$ and $b>0$, expansion in principal minors shows
that $aI+bA$ is a P-matrix for every P-matrix $A$.
The first condition implies the second because $\TP_n\subseteq\mathcal P_n$.
\end{proof}

\section{Positive constant terms}\label{sec:constants}

\begin{proof}[Proof of Theorem~\ref{thm:constant}]
If $n\le3$ and $\mu\ne0$, Theorem~\ref{thm:functions} gives
preservation for every $a\ge0$. If $\mu=0$, any $a>0$ works.
For $n\ge4$, a measure of mass $C$ supported on $[\ell,u]$ satisfies
\[
 \sqrt{M_{r_n}M_{-r_n}}\le C(u/\ell)^{r_n/2}.
\]
Theorem~\ref{thm:repair} gives the assertion with
$c=c_n(u/\ell)^{r_n/2}$.
\end{proof}

For instance, for every fixed order $n$, the function
$a+z/(z+1)+z/(z+2)$ preserves P-matrices when $a$ is sufficiently large.
It has two distinct poles, so it is neither affine nor a constant
shift of a single kernel.

\begin{proposition}\label{prop:unrepairable}
Let $0<p<1/2$ and $a\ge0$. The complete Bernstein function
\[
 f(z)=a+\frac{z^p}{1+z^p}
\]
does not send every strictly totally positive matrix of order $n$ to a
P-matrix for any $n\ge4$. In order four, the family \eqref{eq:Ar} satisfies
\begin{equation}\label{eq:unrepairable}
 r^{-(1-2p)}\det f(A_r)[14,14]
 \longrightarrow-\frac{\alpha^p-\beta^p}{2\sqrt5}<0.
\end{equation}
Thus adding any fixed nonnegative constant fails to restore P-preservation,
even though $f$ is bounded on the positive axis.
\end{proposition}

\begin{proof}
Put $g(z)=z^p/(1+z^p)$. On the principal branch,
\begin{equation}\label{eq:powerkernelmeasure}
 g(z)=\int_0^\infty\frac{z}{z+s}\,
 \frac{\sin(\pi p)}{\pi}
 \frac{s^{p-1}}{1+2s^p\cos(\pi p)+s^{2p}}\,ds.
\end{equation}
Apply Cauchy's formula to
$g(z)/z=z^{p-1}/(1+z^p)$ around the negative real axis. Its upper boundary
value at $-s$ has imaginary part
\[
 -\frac{s^{p-1}\sin(\pi p)}
        {1+2s^p\cos(\pi p)+s^{2p}}.
\]
The small and large circular arcs vanish, giving
\eqref{eq:powerkernelmeasure}. Its measure is positive and has total
mass $1$, since $g(x)\to1$ as $x\to\infty$. Hence $g$ belongs to
\eqref{eq:kernel} and $f=a+g$ is a complete Bernstein function.

The spectral estimates \eqref{eq:spec} give
\[
 g(\lambda_1)\longrightarrow1,\qquad
 g(\lambda_2)\longrightarrow\frac12,\qquad
 g(\lambda_3)-g(\lambda_4)
  =(\alpha^p-\beta^p)r^{-2p}(1+o(1)).
\]
Using \eqref{eq:w} and \eqref{eq:wsums} to group the four middle terms
in \eqref{eq:functionminor}, we obtain
\begin{align*}
 \det g(A_r)[14,14]
 ={}&-\frac r{\sqrt5}
       \bigl(g(\lambda_1)-g(\lambda_2)\bigr)
       \bigl(g(\lambda_3)-g(\lambda_4)\bigr)\\
 &+O(r^{-2p}+r^{-1}).
\end{align*}
Consequently the limit in \eqref{eq:unrepairable} holds first for $a=0$.

The diagonal entries of $A_r(A_r+sI)^{-1}$ belong to $(0,1)$, by
\eqref{eq:res1} and \eqref{eq:res2}. The mass-one representation
\eqref{eq:powerkernelmeasure} therefore gives
$0<g(A_r)_{ii}<1$. Thus
\[
 \det f(A_r)[14,14]
 =a^2+a\tr g(A_r)[14,14]+\det g(A_r)[14,14],
\]
and the first two terms are bounded independently of $r$.
Since $1-2p>0$, they disappear after multiplication by $r^{-(1-2p)}$.
This proves \eqref{eq:unrepairable} for every $a\ge0$.
For each fixed sufficiently large $r$, the higher-dimensional counterexamples
follow by adjoining an identity block and approximating the resulting
nonsingular TN matrix by strictly TP matrices, using continuity of the
matrix function.
\end{proof}

\section{Necessary moments}\label{sec:necessary}

We first construct P-matrices with two long Jordan blocks.

\begin{lemma}\label{lem:jordan}
Let $m\ge2$, let $N$ be the $m\times m$ matrix with ones on its first
superdiagonal and zeros elsewhere, and let $r>0$. Put
\[
 S=r(I-rN)^{-1},\qquad
 T=r^{-1}(I-rN^{\mathsf T})=S^{-\mathsf T},
 \qquad
 J=\begin{pmatrix}1&1\\1&1\end{pmatrix},\quad
 H=\begin{pmatrix}1&-1\\-1&1\end{pmatrix}.
\]
Then
\begin{equation}\label{eq:jordan}
 B_{m,r}=\frac{S\otimes J+T\otimes H}{2}
\end{equation}
belongs to $\mathcal P_{2m}$ and has determinant one. For every
function $h$ holomorphic near $\{r,r^{-1}\}$,
\begin{equation}\label{eq:jfun}
 h(B_{m,r})=\frac{h(S)\otimes J+h(T)\otimes H}{2}.
\end{equation}
\end{lemma}

\begin{proof}
The coordinates in \eqref{eq:jordan} are arranged in $m$ pairs.
Applying the orthogonal change of coordinates
\[
 \frac1{\sqrt2}\begin{pmatrix}1&1\\1&-1\end{pmatrix}
\]
in every pair, and then collecting the first and second coordinates,
gives $S\oplus T$. This proves \eqref{eq:jfun}. It also shows that
the eigenvalues of $B_{m,r}$ are $r$ and $r^{-1}$, each with
multiplicity $m$, and that $\det B_{m,r}=1$.

Fix a principal submatrix of $B_{m,r}$. A pair contributes zero, one,
or two selected coordinates. Divide the represented pairs into
maximal consecutive intervals, separated by pairs with no selected
coordinate. Since $S$ is upper triangular and $T$ has only a diagonal
and a first subdiagonal, the selected matrix is block upper triangular
across these intervals. Its determinant is the product of the
corresponding interval determinants. We may therefore work in one
interval with no omitted pair. The restrictions of $S$ and $T$ to
such an interval have the same forms with a smaller shift matrix;
we retain the notation $S,T,N$ for these restrictions.

Let $F$ be the set of pairs with two selected coordinates and $U$ the
set with one selected coordinate. Assign $\varepsilon_i=1$ or $-1$
to a singleton pair according as its first or second coordinate is
selected, and put $E=\diag(\varepsilon_i:i\in U)$. On the full pairs,
the orthogonal change of coordinates above makes the block to be
eliminated equal to $S[F,F]\oplus T[F,F]$. Its determinant is
$r^{|F|}r^{-|F|}=1$. The Schur complement on the singleton coordinates
is
\begin{equation}\label{eq:jschur}
 \frac{S/F+E(T/F)E}{2},
\end{equation}
where $/F$ denotes the Schur complement with respect to $F$.
If $U$ is empty, the principal determinant is already one.

The inverse block identity gives
\[
 (S/F)^{-1}=S^{-1}[U,U]=r^{-1}(I-rN)[U,U].
\]
Hence $S/F$ is block diagonal on the maximal consecutive runs of
singleton pairs. On a run of length $\ell$, its block is
$r(I-rN_\ell)^{-1}$. Eliminating indices from the lower bidiagonal
matrix $T$ produces only lower triangular connections between these
runs: each new entry joins singleton pairs along a downward path through
eliminated full pairs, so it cannot change a diagonal run block.
Its diagonal block on a run is
$r^{-1}(I-rN_\ell^{\mathsf T})$. Thus \eqref{eq:jschur} is block lower
triangular on the runs. It remains to show that
\[
 W_\ell=\frac{r(I-rN_\ell)^{-1}
       +r^{-1}E(I-rN_\ell^{\mathsf T})E}{2}
\]
has positive determinant for every choice of the singleton signs.

Set $R=r^2$ and $\delta_i=\varepsilon_i\varepsilon_{i+1}$.
Multiplication by $2r(I-rN_\ell)$ gives a tridiagonal matrix $Z_\ell$
whose diagonal entries are $R+1+R\delta_i$ for $i<\ell$ and $R+1$
for $i=\ell$, and whose upper and lower off-diagonal entries are
$-r$ and $-r\delta_i$, respectively. Moreover,
\[
 \det W_\ell=(2r)^{-\ell}\det Z_\ell.
\]
Write $D_i$ for the trailing principal determinant of $Z_\ell$ on
$i,\ldots,\ell$, with $D_{\ell+1}=1$. Then $D_\ell=R+1$ and
\[
 D_i=(R+1)D_{i+1}+R\delta_i(D_{i+1}-D_{i+2}).
\]
Inductively, $D_i>D_{i+1}>0$. Indeed, if $\delta_i=1$, then
\[
 D_i-D_{i+1}=R D_{i+1}+R(D_{i+1}-D_{i+2})>0;
\]
if $\delta_i=-1$, the same difference equals $R D_{i+2}>0$.
The initial difference is $D_\ell-D_{\ell+1}=R>0$.
Consequently every run determinant is positive. Reversing the Schur
complements and the interval decomposition proves that every
principal minor of $B_{m,r}$ is positive.
\end{proof}

Write a complete Bernstein function vanishing at zero as
\begin{equation}\label{eq:zg}
 g(z)=bz+\int_{(0,\infty)}\frac{z}{z+s}\,\mu(ds),
 \qquad b\ge0,\qquad
 \int_{(0,\infty)}\frac{\mu(ds)}{1+s}<\infty.
\end{equation}
Put $C=M_0=\mu((0,\infty))$, allowing $C=\infty$.

\begin{proposition}\label{prop:necessary}
Let $n\ge4$, $a\ge0$, and let $g$ be as in \eqref{eq:zg}, with
$a+g$ nonzero. Suppose that $a+g$ preserves every matrix in $\Pn$.
Then either $\mu=0$, or $b=0$, $C<\infty$, and
\begin{equation}\label{eq:momentlb}
 M_k,M_{-k}<\infty,\qquad
 \sqrt{M_kM_{-k}}\le2a+C
 \quad\left(1\le k\le\left\lfloor\frac n2\right\rfloor-1\right).
\end{equation}
In particular, every unbounded complete Bernstein function preserving
$\Pn$ is affine.
\end{proposition}

\begin{proof}
Suppose that $\mu\ne0$.

When $n>4$, append an identity block to the four-dimensional tests
below; the displayed output minor is unchanged. Take $m=2$ in
Lemma~\ref{lem:jordan}. Since
$S=rI+r^2N$ and $T=r^{-1}I-N^{\mathsf T}$, the principal minor on
coordinates $1,4$ is
\begin{equation}\label{eq:jtest}
 \det(aI+g(B_{2,r}))[\{1,4\},\{1,4\}]
 =\left(a+\frac{g(r)+g(r^{-1})}{2}\right)^2
       -\frac{r^2}{4}g'(r)g'(r^{-1}).
\end{equation}
All limits below are taken as $r\to\infty$.

If $b>0$, dominated and monotone convergence give
\[
 \frac{g(r)}r\longrightarrow b,\qquad
 g'(r)\longrightarrow b,\qquad
 g(r^{-1})\longrightarrow0,\qquad
 g'(r^{-1})\longrightarrow b+M_{-1}>b.
\]
The last limit is allowed to be infinite. For the first two limits,
one may use, for $r\ge1$,
\[
 \frac1{r+s}\le\frac1{1+s},\qquad
 \frac{s}{(r+s)^2}\le\frac1{1+s}.
\]
For $g(r^{-1})\to0$, use
$x/(x+s)\le(1+s)^{-1}$ when $0<x\le1$.
After division by $r^2$, the right side of \eqref{eq:jtest} tends to
$-bM_{-1}/4<0$, with the value $-\infty$ permitted. This contradicts
the assumed preservation.

Suppose next that $b=0$ and $C=\infty$. Then $g(r)\to\infty$.
For $T\ge1$, put
\[
 E_T=\int_{(T,\infty)}s^{-1}\,\mu(ds).
\]
The measure condition in \eqref{eq:zg} gives $E_T<\infty$ and
$E_T\to0$. Splitting the integral at $T$ and applying the
Cauchy--Schwarz inequality gives
\begin{align*}
 g(r)
 &\le\mu((0,T])
   +\left(\int_{(T,\infty)}\frac{r^2s}{(r+s)^2}\,\mu(ds)\right)^{1/2}
      E_T^{1/2}\\
 &\le\mu((0,T])+\bigl(r^2g'(r)\bigr)^{1/2}E_T^{1/2}.
\end{align*}
The first term is finite, and $E_T>0$ because $C=\infty$. Dividing
by $g(r)$ and letting $r\to\infty$ gives
\[
 \liminf_{r\to\infty}\frac{r^2g'(r)}{g(r)^2}\ge E_T^{-1}.
\]
Letting $T\to\infty$, we obtain
\begin{equation}\label{eq:tailratio}
 \frac{r^2g'(r)}{g(r)^2}\longrightarrow\infty.
\end{equation}
Also $g(r^{-1})\to0$ by the same domination used above, and
$g'(r^{-1})\ge g'(1)>0$ for $r\ge1$. Thus division of
\eqref{eq:jtest} by $g(r)^2$ gives a limit of $-\infty$, again a
contradiction.

It follows that $b=0$ and $0<C<\infty$. Fix an integer $k\ge1$ with
$2k+2\le n$, and apply Lemma~\ref{lem:jordan} with $m=k+1$.
Finite geometric expansion of the resolvents gives
\[
 K_s(S)=\sum_{j=0}^{k}\frac{r^{j+1}s^j}{(r+s)^{j+1}}N^j,
 \qquad K_s(z)=\frac{z}{z+s},
\]
whereas the coefficient of $(N^{\mathsf T})^j$ in $K_s(T)$ is
$-s/(r^{-1}+s)^{j+1}$ for $j\ge1$. Consequently, if
\[
 U_k(r)=\int_{(0,\infty)}
       \left(\frac r{r+s}\right)^{k+1}s^k\,\mu(ds),\qquad
 V_k(r)=\int_{(0,\infty)}
       \frac{s}{(r^{-1}+s)^{k+1}}\,\mu(ds),
\]
then $g(S)_{1m}=U_k(r)$ and $g(T)_{m1}=-V_k(r)$.
Both integrals are finite for each fixed $r$ because $\mu$ is finite
and the integrands are bounded. The diagonal entries of $g(B_{m,r})$
are all $(g(r)+g(r^{-1}))/2$. Formula~\eqref{eq:jfun} therefore yields
\begin{equation}\label{eq:ktest}
 \det(aI+g(B_{m,r}))[\{1,2m\},\{1,2m\}]
 =\left(a+\frac{g(r)+g(r^{-1})}{2}\right)^2
       -\frac14 U_k(r)V_k(r).
\end{equation}

The two integrands increase pointwise with $r$, so monotone
convergence gives
\[
 U_k(r)\longrightarrow M_k,\qquad
 V_k(r)\longrightarrow M_{-k}.
\]
Also $g(r)\to C$ and $g(r^{-1})\to0$. Both moments are strictly
positive since $\mu\ne0$. If either were infinite, the right side
of \eqref{eq:ktest} would tend to $-\infty$. Otherwise its limit is
\[
 \left(a+\frac C2\right)^2-\frac14M_kM_{-k}.
\]
Append an identity block when $2m<n$. Preservation requires this
limit to be nonnegative, proving \eqref{eq:momentlb}.
\end{proof}

\section{Cycle products and polynomial estimates}\label{sec:sufficient}

We now prove Theorem~\ref{thm:kernel}. Fix $n\ge4$ and write
\[
 m=\lfloor n/2\rfloor,\qquad r=m-1=r_n,\qquad
 w(s)=\max(s,s^{-1})^r,\qquad
 K_s(A)=A(A+sI)^{-1}.
\]
For a fixed $A\in\Pn$, write $K_s=K_s(A)$ and use
the polynomial $p_A$ from \eqref{eq:p}.
Every coefficient of $p_A$ is strictly positive.
The same-parameter cycle bound below is the uniform bounded-set form
of Lin and Sturmfels~\cite[Corollary~8]{LS}.
We give a direct proof because its comparison of Hamilton cycles
will also be used with independent kernel parameters.

We regard a permutation \(\sigma\) as mapping a column index \(c\)
to its row index \(\sigma(c)\). Thus a directed edge \(c\to d\)
of a permutation carries the matrix entry in row \(d\), column \(c\).
Cycles are identified up to cyclic rotation.

\begin{lemma}\label{lem:graph}
Fix a finite vertex set \(I\) of size \(j\), with nonnegative weights
on its directed edges. Suppose that the products along all proper
simple cycles are at most \(B\), where \(B\ge1\).
Let \(\gamma\) maximize the product \(M\) among Hamilton cycles on
\(I\). If \(M\ge1\), the product along every other Hamilton cycle is
bounded by a constant depending only on \(j\) and \(B\).
\end{lemma}

\begin{proof}
Let \(\eta\ne\gamma\). If its product is zero there is nothing to prove.
Choose an edge \(v\to w\) of \(\eta\) which is not an edge of \(\gamma\).
Together with the path from \(w\) to \(v\) along \(\gamma\), it forms
a proper simple cycle, of length between two and \(j-1\).
Delete these edges, with multiplicity, from the union of
\(\gamma\) and \(\eta\). The remaining directed multigraph is balanced
and has fewer than \(2j\) edges. It can be decomposed into simple
directed cycles. At most one of them has length \(j\), and its product
is at most \(M\). Every other cycle product is at most \(B\).
Multiplying all edge weights gives
\[
 M\,\prod_{e\in\eta}h_e\le B^{2j}M.
\]
Division by \(M\) proves the assertion.
\end{proof}

\begin{lemma}\label{lem:samecycle}
For every \(d\), there is a constant \(b_d\) with the following
property. If every principal minor of a real \(d\times d\) matrix
has absolute value at most one, then every simple directed cycle
product of its entries has absolute value at most \(b_d\).
\end{lemma}

\begin{proof}
Induct on the length \(j\) of the cycle. Diagonal entries give the
case \(j=1\). Fix a set \(I\) of size \(j\) and suppose all shorter
cycle products are bounded. In the expansion of \(\det X[I,I]\),
the terms whose permutations have only proper cycles are bounded.
Let \(M\) be the largest absolute Hamilton product. If \(M<1\),
there is nothing further to prove. Otherwise
Lemma~\ref{lem:graph} bounds every other Hamilton product.
Since \(\abs{\det X[I,I]}\le1\), the remaining maximal Hamilton term
is bounded as well. There are only finitely many vertex sets and
permutations in the fixed dimension.
\end{proof}

We shall repeatedly use a factor-specific polynomial norm estimate.
For fixed degree bounds and a fixed number of real polynomials,
\begin{equation}\label{eq:productnorm}
 \prod_{\nu=1}^k\max_h\abs{[z^h]P_\nu(z)}
 \le C\sup_{1\le z\le2}\abs{\prod_{\nu=1}^kP_\nu(z)}.
\end{equation}
Indeed, normalize every coefficient maximum to one. The resulting
coefficient unit spheres are compact, and their product is never the
zero polynomial. Its supremum on \([1,2]\) therefore has a positive
minimum. This proves \eqref{eq:productnorm}.
Ordinary polynomial coefficient norm equivalence follows in the
same way.

Let \(q(s)=\sum_{h=0}^na_hs^h\), with every \(a_h>0\).
On each interval \([x,2x]\), the functions \(q(s)w(s)^p\) vary by
a bounded factor when \(1\le p\le n\). Consequently, for polynomials
of degree at most \(n\), every coefficient projection is uniformly
bounded in the norm
\begin{equation}\label{eq:weightednorm}
 \sup_{s>0}\frac{\abs{P(s)}}{q(s)w(s)^p}.
\end{equation}
For example, coefficient norm equivalence on \([x,2x]\) gives
\[
 \abs{c_h}x^h
 \le C_n\norm{P/(qw^p)}_\infty q(x)w(x)^p.
\]
Taking all \(x>0\), and then summing any prescribed subset of
coefficients, proves the assertion. The same bound applies to a
coordinate projection of a polynomial in several variables, with
all other variables held fixed.

\begin{lemma}\label{lem:capacity}
Let \(h_1,\ldots,h_k\in\{1,\ldots,n-1\}\), let
\(D=\sum h_i\), and write
\[
 \alpha=D/k=b+\rho/k,\qquad b=\lfloor\alpha\rfloor,\quad
 0\le\rho<k.
\]
Put \(\theta_D=\sup_{x>0}x^D/q(x)^k\).
For any \(\epsilon_i\in\{0,1\}\) with \(\sum\epsilon_i=\rho\),
\[
 \frac{\prod_i s_i^{h_i}}{\prod_iq(s_i)}
 \le C_{n,k}\theta_D
       \prod_i s_i^{h_i-b-\epsilon_i}.
\]
\end{lemma}

\begin{proof}
The positive function \(q(x)/x^\alpha\) attains its minimum at
some \(x_0>0\). The numbers \(a_\ell x_0^\ell/q(x_0)\) form a
probability distribution of mean \(\alpha\). Discrete convexity gives
\[
 \frac{q(x_0y)}{q(x_0)}
 \ge y^b(1-\delta+\delta y),\qquad \delta=\rho/k.
\]
If \(\delta=0\), this immediately proves the claim. Otherwise use
\(1-\delta+\delta y\ge1-\delta\) in the slots with
\(\epsilon_i=0\), and
\((1-\delta+\delta y)/y\ge\delta\) in the other slots.
The resulting constant is
\((1-\delta)^{-(k-\rho)}\delta^{-\rho}\), and
\(x_0^D/q(x_0)^k=\theta_D\).
\end{proof}

\section{A weighted scalar product estimate}

\begin{lemma}\label{lem:scalar}
Let \(3\le L\le n\), let \(p_1,\ldots,p_L\) be positive integers,
and let \(f_i\) be real polynomials supported in
\([p_i,n-p_i]\). Set \(V_i=q w^{p_i}\), and put
\[
 \mathcal A=\prod_i\norm{f_i/V_i}_\infty,\qquad
 \Lambda=\sup_{s>0}\frac{\abs{\prod_if_i(s)}}{q(s)^L},
\]
\[
 \mathcal M=\max_{a<b}\ 
 \sup_{s,\ s_i:\ i\ne a,b}
 \frac{\abs{f_a(s)f_b(s)\prod_{i\ne a,b}f_i(s_i)}}
 {V_a(s)V_b(s)\prod_{i\ne a,b}V_i(s_i)}.
\]
Then \(\mathcal A\le C_n(\Lambda+\mathcal M)\).
No restriction on \(\sum p_i\) is required.
\end{lemma}

\begin{proof}
If some factor is zero, the assertion is immediate. Otherwise each
\(p_i\le n/2\), so only finitely many support intervals arise for fixed
\(n\). Suppose the
claim fails and divide every factor by its individual weighted norm.
We obtain a sequence with all individual norms equal to one, while
\(\Lambda,\mathcal M\to0\). In particular, every normalized pair product
has supremum tending to zero.

Write
\[
 \kappa_{p,h}=\sup_{s>0}\frac{s^h}{q(s)w(s)^p}.
\]
The triangle inequality allows us to choose a coefficient \(c_i\)
of degree \(h_i\) in every factor so that
\begin{equation}\label{eq:selected}
 \abs{c_i}\kappa_{p_i,h_i}\ge\frac1{n+1}.
\end{equation}
Pass to a subsequence on which all \(h_i\) are fixed.
Applying \eqref{eq:productnorm} on each dyadic interval gives
\begin{equation}\label{eq:vanishcapacity}
 \abs{\prod_i c_i}\theta_D\longrightarrow0,\qquad
 D=\sum_i h_i,\quad \theta_D=\sup_x x^D/q(x)^L,
\end{equation}
and, for every selected pair,
\begin{equation}\label{eq:vanishpair}
 \abs{c_ac_b}\sup_x
 \frac{x^{h_a+h_b}}{q(x)^2w(x)^{p_a+p_b}}
 \longrightarrow0.
\end{equation}
Here the factor-specific estimate controls products of individual
coefficients; it is not a coefficient extraction from a possibly
cancelling product polynomial.

Call a tuple good if an allocation in Lemma~\ref{lem:capacity}
makes every residual \(h_i-b-\epsilon_i\) lie in \([-r,r]\).
For a good tuple, that lemma gives
\[
 \prod_i\kappa_{p_i,h_i}
 \le\prod_i\kappa_{1,h_i}\le C_n\theta_D,
\]
contradicting \eqref{eq:selected}--\eqref{eq:vanishcapacity}.
Thus the selected tuple is bad. Simultaneous reflection
\[
 q(s)\mapsto s^nq(1/s),\qquad
 f_i(s)\mapsto s^nf_i(1/s),\qquad h_i\mapsto n-h_i
\]
preserves all norms, capacities and support intervals. We may assume
\(\alpha=D/L\le n/2\).

We need two rounding facts: a bad tuple in this range has at least
two indices \(h_i\le m-1\), and there is an allocation with every
residual at least \(-r\).
For even \(n=2m\), equality \(\alpha=m\) makes every tuple good.
Otherwise \(b\le m-1\), so every allocation has all residuals at least \(-r\).
If at most one index is below \(m\), and \(H=\max h_i\), then
\[
 H-\alpha\le
 \frac{(L-1)(2m-1)-(L-2)m-1}{L}=m-1.
\]
All upper bounds then hold before any ceil slots are chosen, and
the tuple is good.

For odd \(n=2m+1\), the residual lower bound \(-r\) is automatic
when \(b\le m-1\).
If \(b=m\), write \(D=Lm+\rho\), where \(2\rho\le L\).
Only indices equal to one must avoid a ceil. If there were more
than \(L-\rho\) of these, then
\[
 D\le L+(\rho-1)(2m-1)<Lm+\rho,
\]
because the difference is at most
\((m-1)(2\rho-L)-(2m-1)<0\).
Thus an allocation with all residuals at least \(-r\) exists.
If at most one index is below \(m\) and \(H\le2m-1\), the preceding
maximum estimate again makes the tuple good. If \(H=2m\), let \(c\)
be the number of such maximal indices. Then $D>Lm$, while $D/L\le m+1/2$, so $b=m$ and
\[
 D\ge Lm+(c-1)m+1,\qquad
 \rho\ge(c-1)m+1\ge c.
\]
Assign a ceil to all \(c\) maximal indices and avoid the possible
one index equal to one; this is possible since
\(\rho\le L/2\le L-1\). Every other upper bound already holds.
This proves both rounding facts for arbitrary \(L\).

Normalize \(q(1)=1\), adjusting the factors by the same scalar,
and put \(\beta=q'(1)\). The logarithmic slope
\(\beta(s)=s q'(s)/q(s)\) is increasing: its derivative with respect
to \(\log s\) is the variance of the distribution proportional to
\(a_ks^k\). The logarithmic derivatives of the selected monomial
ratios below and above one are, respectively,
\begin{equation}\label{eq:slopes}
 h_i+p_ir-\beta(s),\qquad h_i-p_ir-\beta(s).
\end{equation}

If \(n=2m\) and \(\beta\ge m\), every selected monomial norm is
attained on \(s\le1\), since \(h_i-p_ir\le h_i-r\le m\).
Choose an allocation with all residuals at least \(-r\).
At maximizing points \(s_i\le1\),
its residuals \(e_i\ge-r\) give
\[
 s_i^{e_i}/w(s_i)^{p_i}=s_i^{e_i+p_ir}\le1.
\]
Lemma~\ref{lem:capacity} contradicts
\eqref{eq:selected}--\eqref{eq:vanishcapacity}.
If \(\beta<m\), the weighted monomial ratio for each low degree
\(h_i\le m-1\) attains its maximum at one by \eqref{eq:slopes}.
The two low indices therefore contradict
\eqref{eq:selected} and \eqref{eq:vanishpair}.

For odd \(n=2m+1\), if \(\beta\le m\), every selected low degree
again has its weighted monomial ratio maximized at one:
\(h_i+p_ir\ge m\) and \(h_i-p_ir\le0\).
If \(\beta\ge m+1\), every selected monomial ratio attains its maximum
on \(s\le1\), since \(h_i-p_ir\le m+1\).
The preceding pair-overlap and capacity arguments apply, respectively.
It remains to consider \(m<\beta<m+1\).
The weighted monomial ratio for each selected degree
\(2\le h_i\le m-1\) attains its maximum at one.
A selected degree one necessarily has \(p_i=1\), because \(h_i\ge p_i\).
We claim
\begin{equation}\label{eq:massgap}
 \min\left(\kappa_{1,1},
       \max_{2\le h\le2m}\kappa_{1,h}\right)
 \le 2(m+1)(m+2).
\end{equation}
Indeed, put \(P=\sum_{k\le m}a_k\). Since \(\beta\le m+1\),
the mass above \(m+1\) is at most \((m+1)P\).
If \(P\ge1/[2(m+2)]\), one coefficient with \(k\le m\) is at least
\(1/[2(m+1)(m+2)]\). On \(s\le1\) it bounds \(s^m/q(s)\),
and on \(s\ge1\) the degree-one normalized ratio is at most one.
This bounds \(\kappa_{1,1}\). Otherwise \(a_{m+1}>1/2\), and that
term bounds every \(\kappa_{1,h}\), \(h\ge2\), by two.

There are at least two selected low indices. If two equal one,
both have \(p_i=1\), and their identical monomial ratio functions
contradict \eqref{eq:vanishpair}. If two are at least two, their
center overlap gives the same contradiction. In the remaining case
there is exactly one degree-one factor and another low factor.
If \(\kappa_{1,1}\) is bounded along a subsequence, pair these at one.
Otherwise pass to a subsequence on which \(\kappa_{1,1}\to\infty\).
Then \eqref{eq:massgap} bounds all selected weighted monomial
norms of degrees at least two, since
\(\kappa_{p,h}\le\kappa_{1,h}\).
There are at least two such factors because \(L\ge3\).
Their overlap at one contradicts \eqref{eq:vanishpair}.
All cases are impossible.
\end{proof}

\section{Fixed-column determinants}

Fix \(A\in\Pn\), a principal set \(I\) of size \(j\), and an ordered
partition \(\mathcal C=(C_1,\ldots,C_L)\) of \(I\), with
\(p_g=\abs{C_g}\ge1\). Let
\[
 U_{I,\mathcal C}(\boldsymbol s)
 =\det\bigl[(K_{s_g})_{i,c}\bigr]_
       {i\in I,\ c\in C_g},
\]
The columns retain their original order, and the parameter assigned
to each column is fixed by \(\mathcal C\).
The target is
\begin{equation}\label{eq:fixedtarget}
 \abs{U_{I,\mathcal C}(\boldsymbol s)}
 \le C_n\prod_gw(s_g)^{p_g}.
\end{equation}

Since \(A^{-1}\in\Pn\), Jacobi gives
\[
 \det K_s[R,C]
 =\pm\frac{\det(I+sA^{-1})[C^c,R^c]}
               {\det(I+sA^{-1})}.
\]
It follows that
\begin{equation}\label{eq:groupdegree}
 p_A(s)\det K_s[R,C]
 \quad\hbox{is a polynomial of degree at most }n-p
 \quad(\abs R=\abs C=p).
\end{equation}
When \(\ell=\abs{R\setminus C}>0\), all coefficients below degree
\(\ell\) vanish. Every nonempty principal minor of \(K_s\) lies in \((0,1)\):
the Jacobi numerator for a principal minor is a partial sum of the
positive principal-minor expansion of its denominator.

The case \(L=1\) of \eqref{eq:fixedtarget} is therefore immediate.
For \(L=2\), expand by the fixed column set \(C_1\), of size \(p\);
write \(q_0=j-p\). The expansion is
\[
 U=\sum_{\substack{R\subset I\\\abs R=p}}\pm
 \det K_s[R,C_1]\det K_t[I\setminus R,C_2].
\]
The term \(R=C_1\) is a product of principal minors and is bounded by
one. Otherwise \(\ell=\abs{R\setminus C_1}\ge1\), and the two
numerators have degree ranges \([\ell,n-p]\) and
\([\ell,n-q_0]\). At a common parameter their minor product is a
sum of products along disjoint cycles, hence bounded by
Lemma~\ref{lem:samecycle}. Applying \eqref{eq:productnorm} on
\([x,2x]\), every individual coefficient pair satisfies
\[
 \abs{c_hd_v}x^{h+v}\le C_n p_A(x)^2.
\]
Apply Lemma~\ref{lem:capacity}, assigning the ceil to the higher degree
when necessary. The term is bounded by
\[
 C_n\max(s/t,t/s)^{\lfloor\abs{h-v}/2\rfloor}
 \le C_nw(s)w(t).
\]
Here \(\abs{h-v}\le n-2\), and
\(\lfloor(n-2)/2\rfloor=r\).
Summing proves \eqref{eq:fixedtarget} for \(L=2\).
In particular, this estimate holds for each fixed column partition.

Put
\[
 N(\boldsymbol s)=\prod_gp_A(s_g)\,
                         U_{I,\mathcal C}(\boldsymbol s),\qquad
 V_g(s)=p_A(s)w(s)^{p_g}.
\]
Grouped Laplace expansion and \eqref{eq:groupdegree} show that
\(N\) has degree at most \(n-p_g\) in its \(g\)-th variable.
Let \(\Pi_g\) retain its degrees \(p_g,\ldots,n-p_g\), and let
\(\Pi=\prod_g\Pi_g\). Each of these coordinate projections is
uniformly bounded in the norm
\[
 \norm{P}_{\boldsymbol V}
 =\sup_{\boldsymbol s}
       \frac{\abs{P(\boldsymbol s)}}{\prod_gV_g(s_g)}
\]
by the coefficient projection estimate following
\eqref{eq:weightednorm}.

\begin{lemma}\label{lem:lowjets}
Suppose \eqref{eq:fixedtarget} is known in every principal order
strictly below \(j\), for the same source dimension \(n\).
Then
\[
 \norm{N-\Pi N}_{\boldsymbol V}\le C_n.
\]
\end{lemma}

\begin{proof}
Fix a group \(C=C_g\) of size \(p\), and \(0\le h<p\).
For \(B\subset C\), let \(U_B(s)\) be the determinant with \(K_s\)
columns on \(B\), identity columns on \(C\setminus B\), and the
other columns unchanged. Holding the other parameters fixed, write
\(K_s=I+Z_s\), where \(Z_s=O(s)\) as \(s\to0\).
Multilinearity gives
\begin{equation}\label{eq:jet}
 U(s)-\sum_{\substack{B\subset C\\\abs B\le h}}
 (-1)^{h-\abs B}
 \binom{p-\abs B-1}{h-\abs B}\,U_B(s)
 =O(s^{h+1}).
\end{equation}
Indeed, keep at most \(h\) columns \(Z_s\) and expand them as \(K_s-I\).
For a fixed subset \(B\) of size \(b\le h\), its coefficient is
\[
 \sum_{d=0}^{h-b}(-1)^d\binom{p-b}{d}
 =(-1)^{h-b}\binom{p-b-1}{h-b}.
\]
An identity column indexed by \(c\) has its only nonzero entry in
row \(c\). Deleting these rows and columns leaves a principal
fixed-column determinant on \(I\setminus(C\setminus B)\), of order
\(j-p+\abs B<j\), still formed from the kernels \(K_s(A)\).
The induction hypothesis gives a factor \(w(s)^{\abs B}\),
which is at most \(w(s)^p\) because \(w\ge1\).
Multiply \eqref{eq:jet} by the original product of \(p_A\)'s.
Every term is a polynomial of fixed degree, including when \(B\)
is empty and the unused \(p_A(s)\) is simply an extra factor.
Coefficient projection therefore bounds every degree-\(h\) slice
of \(N\), \(h<p\), in the full norm \(\norm{\cdot}_{\boldsymbol V}\).
Thus \(\norm{(I-\Pi_g)N}_{\boldsymbol V}\le C_n\).
The telescoping identity for \(I-\prod_g\Pi_g\), and boundedness of
the other coordinate projections, prove the claim.
\end{proof}

\section{Double induction and simultaneous diagonal scaling}

\begin{proposition}\label{prop:fixedall}
The estimate \eqref{eq:fixedtarget} holds for every \(1\le j\le n\)
and every fixed column partition. Moreover, the products
\[
 \prod_{c\in I}H_{\sigma(c),c},\qquad
 H_{ab}=\sup_{s>0}\frac{\abs{(K_s)_{ab}}}{w(s)},
\]
are uniformly bounded for all simple directed cycles \(\sigma\).
\end{proposition}

\begin{proof}
For each fixed \(A\), the kernel \(K_s\) extends continuously to
\(s=0\) and \(s=\infty\), with limits \(I\) and \(0\).
Thus every \(H_{ab}\) is finite.
We induct first on \(j\). The diagonal bound gives the case \(j=1\).
Assume the fixed-column bounds and the cycle-product bounds for \(H\)
in every smaller principal order. For this \(j\), induct on \(L\). The cases \(L=1,2\)
were proved in the preceding section.

Suppose that the bound fails for some \(L\ge3\).
There are finitely many index sets and column partitions, so fix
one of each and choose a sequence of matrices along which the
normalized determinant norm tends to infinity. Terms containing
only shorter cycles are bounded by the induction on \(j\).
Consequently, the maximum of
\[
 \prod_{c\in I}H_{\gamma(c),c}
\]
over Hamilton cycles is unbounded. Pass to a subsequence on which
one fixed cycle \(\gamma\) attains this maximum and the maximum tends
to infinity. Lemma~\ref{lem:graph} bounds the products associated
with every other Hamilton cycle.

In the grouped Laplace expansion, exactly one row partition can
contain \(\gamma\):
\[
 R_g=\gamma(C_g)\qquad(1\le g\le L).
\]
Let
\[
 G(\boldsymbol s)=\prod_g\det K_{s_g}[R_g,C_g],
 \qquad f_g(s)=p_A(s)\det K_s[R_g,C_g].
\]
Expanding this product of minors, the permutation \(\gamma\) occurs
once. Every other global permutation is either a proper cycle cover
or a different Hamilton orientation. The corresponding product of entries of \(H\) is bounded.
The same holds for every other row partition.
Consequently, for one sign \(\epsilon\),
\begin{equation}\label{eq:isolation}
 \norm{N-\epsilon\prod_gf_g}_{\boldsymbol V}\le C_n.
\end{equation}
Apply \(\Pi\) to \eqref{eq:isolation}. Since \(\Pi_g\) acts only on
\(s_g\),
\[
 \Pi\left(\prod_gf_g(s_g)\right)
 =\prod_g\widetilde f_g(s_g),\qquad
 \widetilde f_g=\Pi_gf_g.
\]
By Lemma~\ref{lem:lowjets}, \(\Pi N\) is still unbounded in full
normalized norm. Hence the product of the individual weighted norms
of \(\widetilde f_g\) is unbounded. If some \(p_g>n/2\), its retained
degree interval is empty, making this product zero, a contradiction.
Otherwise each factor has support \([p_g,n-p_g]\).

Set \(s_a=s_b=s\). The original determinant then has the column
groups \(C_a\) and \(C_b\) merged, so the induction hypothesis for
\(L-1\) groups bounds \(N/\prod_gV_g\) on this set.
Both \(N-\Pi N\) and
\(\Pi(N-\epsilon\prod_gf_g)\) have already been bounded for
independent parameters. Their restrictions therefore show that
\(\prod_g\widetilde f_g\) has the required bound on the same set.
This bounds \(\mathcal M\) in Lemma~\ref{lem:scalar} for every pair
\(a<b\), with the factors \(\widetilde f_g\).

When all parameters equal \(x>0\), expanding
\(G(x,\ldots,x)\) gives products of entries of \(K_x\) along disjoint
cycles. Lemma~\ref{lem:samecycle} gives
\[
 \abs{\prod_gf_g(x)}\le C_np_A(x)^L.
\]
On \([x,2x]\), \eqref{eq:productnorm} therefore bounds every product
of individual factor coefficients, multiplied by its corresponding
power of \(x\), by \(C_np_A(x)^L\).
Truncation retains a subset of these coefficients. Summing its
finitely many products proves
\[
 \sup_{x>0}\frac{\abs{\prod_g\widetilde f_g(x)}}{p_A(x)^L}\le C_n.
\]
Thus \(\Lambda\) in Lemma~\ref{lem:scalar} is bounded as well.
Applying that lemma with \(q=p_A\) bounds the product of the individual
weighted norms, contradicting its divergence. This proves the fixed-column
bound for every \(L\) in order \(j\).

To complete the induction on \(j\), suppose that a product of entries
of \(H\) along a \(j\)-cycle were unbounded. Choose a maximizing cycle
as above. Every other permutation term is uniformly bounded.
Apply the fixed-column estimate with all \(j\) columns carrying
independent parameters. It bounds the maximal cycle term at every
such parameter tuple. Taking independent suprema bounds the product
of the entries of \(H\) along this cycle, a contradiction.
The simultaneous induction is complete.
\end{proof}

\begin{proof}[Proof of Theorem~\ref{thm:kernel}]
By Proposition~\ref{prop:fixedall}, all simple-cycle products of the
nonnegative matrix \(H=(H_{ij})\) are bounded.
We apply the classical cycle criterion for positive diagonal scaling;
see Hershkowitz and Schneider~\cite[Theorem~2.5]{HS}.
The following path construction gives the required scaling directly.
Choose \(B_n\ge1\) so that the product of the normalized weights
\(H_{ij}/B_n\) around every simple cycle is at most one.
Regard \(j\to i\) as an edge of weight \(H_{ij}/B_n\).
Let \(q_i\) be the maximum of the products along all simple directed
paths ending at \(i\), including the empty path of product one.
These maxima are finite and positive. Appending any edge \(j\to i\)
to a maximizing path for \(q_j\), and removing cycles whose products
are at most one, gives
\[
 q_i\ge (H_{ij}/B_n)q_j.
\]
For \(Q_A=\diag(q_1,\ldots,q_n)\) this yields
\[
 \abs{(Q_A^{-1}K_sQ_A)_{ij}}
 \le H_{ij}w(s)\frac{q_j}{q_i}\le B_nw(s).
\]
The construction depends only on \(A\).
\end{proof}

\section{Positive kernel measures and a uniform scalar repair}

\begin{proof}[Proof of Theorem~\ref{thm:repair}]
For $\mu\ne0$, necessity and the lower bound follow from
Proposition~\ref{prop:necessary} with $k=r_n$.
For sufficiency, suppose that \(W=\int w(s)\,\mu(ds)<\infty\).
The matrix \(Q_A\) from Theorem~\ref{thm:kernel} is independent of \(s\),
so every entry of \(Q_A^{-1}g(A)Q_A\) has absolute value at most
\(B_nW\).
If \(a>nB_nW\), the matrix
\[
 aI+Q_A^{-1}g(A)Q_A
\]
has positive diagonal and is strictly diagonally dominant by rows.
Every principal submatrix has the same property and positive
determinant. To see the last assertion directly, continuously multiply
its off-diagonal part by \(t\in[0,1]\). Strict row dominance prevents
a nonzero null vector, by considering a coordinate of maximal
absolute value. The determinant never vanishes and at \(t=0\) is the
positive product of its diagonal entries.
Positive diagonal similarity preserves principal minors. Thus
\(aI+g(A)\in\Pn\), uniformly in \(A\).

For a nonzero measure with finite \(M_r,M_{-r}\), put
\[
 \tau=(M_r/M_{-r})^{1/(2r)},\qquad
 R=\sqrt{M_rM_{-r}}.
\]
Push \(\mu\) forward under \(s\mapsto s/\tau\). Its two endpoint
moments both equal \(R\), and its weighted mass is at most \(2R\).
The identity \(K_s(\tau A)=K_{s/\tau}(A)\), and the invariance of
\(\Pn\) under positive scalar multiplication, transfer the
preceding repair back to \(\mu\). Thus \(a>2nB_nR\) is sufficient.
The zero measure is immediate by choosing any \(a>0\).
\end{proof}

\begin{proof}[Proof of Proposition~\ref{prop:outerlayers}]
The third condition implies the first.
The proof of Proposition~\ref{prop:necessary} uses only principal
minors of order two, and therefore shows that the first condition
implies the second. If $\mu=0$, the nonzero function $g=bz$ preserves
P-matrices. In the other case, the two finite moments imply
$M_0\le\sqrt{M_{r_n}M_{-r_n}}<\infty$, and
Theorem~\ref{thm:repair} proves the third condition.
\end{proof}

\section{Preservation in every order}\label{sec:orders}

\begin{proof}[Proof of Corollary~\ref{cor:dimensions}]
Suppose that such constants exist and that $\mu\ne0$.
Proposition~\ref{prop:necessary} in order four gives $b=0$ and
$C<\infty$. Applying the same proposition in order $2j+2$ gives
$M_j,M_{-j}<\infty$ for every $j\ge1$.

Conversely, if $\mu=0$, positive scalar shifts of $bA$ are P-matrices
by the principal-minor expansion. If $b=0$ and all the displayed
moments are finite, Theorem~\ref{thm:repair} supplies a finite repair
constant in every order $n\ge4$. Orders at most three follow from
Theorem~\ref{thm:functions}. In the case $g=0$, choose $a_n>0$.
\end{proof}

\begin{proof}[Proof of Corollary~\ref{cor:universal}]
Write $f=a+g$, where $g$ has the form \eqref{eq:zg}. By
Proposition~\ref{prop:necessary}, either $\mu=0$, giving the affine
family, or $b=0$ and $0<C<\infty$. In the latter case, the same
constant $a$ must satisfy
\begin{equation}\label{eq:allmoments}
 \sqrt{M_jM_{-j}}\le2a+C\qquad(j\ge1).
\end{equation}
If $\mu$ is not supported at one point, choose two disjoint compact
intervals $[\ell_1,u_1]$ and $[\ell_2,u_2]$ of positive masses
$c_1,c_2$, with $u_1<\ell_2$. Then
\[
 M_jM_{-j}\ge c_1c_2\left(\frac{\ell_2}{u_1}\right)^j
       \longrightarrow\infty,
\]
contradicting \eqref{eq:allmoments}. Thus $\mu=c\delta_s$ for some
$c,s>0$, which gives the second family.

\begin{samepage}
For the converse, the principal-minor expansion shows that $aI+bA$
is a P-matrix whenever $A$ is a P-matrix and $a,b\ge0$ are not both
zero. Also
\[
 A(A+sI)^{-1}=(I+sA^{-1})^{-1}
\]
is a P-matrix, by closure under inversion and nonnegative scalar
shifts. Positive multiplication and addition of $aI$ preserve this
property. Hence both displayed families preserve every order.
\end{samepage}
\end{proof}

\end{document}